\documentclass[11pt]{article}

\usepackage{times,amssymb,amsmath,exscale,array,latexsym}
\usepackage{graphicx}
\usepackage{epsfig}
\usepackage{color}
\definecolor{marin}{rgb}   {0.,   0.3,   0.7}
\definecolor{rouge}{rgb}   {0.8,   0.,   0.}
\definecolor{sepia}{rgb}   {0.8,   0.5,   0.}
\usepackage[colorlinks,citecolor=marin,linkcolor=rouge,
            bookmarksopen,
            bookmarksnumbered
           ]{hyperref}

\newcommand{\e}{\ensuremath{\mathrm{e}}}

\newtheorem{lemma}{Lemma}[section]

\newtheorem{proposition}[lemma]{Proposition}

\numberwithin{equation}{section}

\newcommand{\QED}{\mbox{}\hfill \raisebox{-0.2pt}{\rule{5.6pt}{6pt}\rule{0pt}{0pt}}
          \medskip\par}
\newenvironment{Proof}{\noindent
    \parindent=0pt\abovedisplayskip = 0.5\abovedisplayskip
    \belowdisplayskip=\abovedisplayskip{\bfseries Proof. }}{\QED}

\title{Applying splitting methods with complex coefficients to the numerical integration of unitary problems} 

\author{S. Blanes$^{1}$, F. Casas$^{2}$, A. Escorihuela-Tom\`as$^{3}$ \\[2ex]
$^{1}$ {\small\it Universitat Polit\`ecnica de Val\`encia, Instituto de Matem\'atica Multidisciplinar, 46022-Valencia, Spain}\\{
\small\it email: serblaza@imm.upv.es}\\[1ex]
$^{2}$ {\small\it Departament de Matem\`atiques and IMAC, Universitat Jaume I, 12071-Castell\'on, Spain}\\{
\small\it email: Fernando.Casas@mat.uji.es}\\[1ex]
$^{3}$ {\small\it Departament de Matem\`atiques, Universitat Jaume I, 12071-Castell\'on, Spain}\\{
\small\it email: alescori@uji.es}\\[1ex]
}

\begin{document}
\mathsurround 0.8mm

\maketitle

\begin{abstract}

We explore the applicability of splitting methods involving complex coefficients to solve numerically 
the time-dependent Schr\"odinger equation. We prove that a particular class of integrators are
conjugate to unitary methods for sufficiently small step sizes when applied to problems defined in the group $\mathrm{SU}(2)$. In the general case,
the error in both the energy and the norm of the numerical approximation provided by these methods
does not possess a secular component over long time intervals, when
combined with pseudo-spectral discretization techniques in space.

\end{abstract}


\section{Introduction}

Splitting methods constitute a natural choice for the numerical time integration of differential equations
of the form
\begin{equation} \label{gen.eq}
   \frac{du}{dt} = A(u) + B(u),  \qquad u(0) = u_0,
\end{equation}   
when each subproblem 
\[
    \frac{du}{dt} = A(u), \qquad  \frac{du}{dt} = B(u)
\]
with $u(0)=u_0$ can be solved explicitly \cite{blanes16aci,hairer06gni,mclachlan02sm}. Then, by composing the solution of each part with appropriately chosen
coefficients, it is possible to construct an integrator of a given order $r \ge 1$ for (\ref{gen.eq}). In the particular case
of a linear problem,
\begin{equation} \label{evol.1}
  \frac{du }{dt} = A u + B u, 
\end{equation}
a splitting method is a composition of the form
\begin{equation}   \label{td.2b}
 \Psi^{[r]}(h) = \e^{b_{s+1} h B}  \e^{a_s h A}  \e^{b_s h B} \cdots  \e^{a_1 h A}  \, \e^{b_1 h B},
\end{equation}
where $h := \Delta t$ is the time step and the coefficients $a_j$, $b_j$ are 
chosen as solutions of the order conditions, a set of polynomial equations that must be satisfied to achieve
an order of accuracy $r$, i.e., so that $\exp(h(A+B)) u_0 - \Psi^{[r]}(h) u_0 = \mathcal{O}(h^{r+1})$. 

The simplest
example within this class is the Lie--Trotter splitting,
\begin{equation} \label{lie-trotter}
   \e^{h A} \, \e^{h B} \qquad \mbox{ or } \qquad \e^{h B} \, \e^{h A},
\end{equation}
providing a first order approximation ($r=1$), whereas the palindromic versions
\begin{equation} \label{strang}
 \mathcal{S}(h) = \e^{h/2 A} \, \e^{h B} \, \e^{h/2 A}   \qquad \mbox{ or } \qquad \mathcal{S}(h) = \e^{h/2 B} \, \e^{h A} \, \e^{h/2 B},
\end{equation}
known as  Strang splittings, are methods of order $r=2$. 

Although very efficient high order splitting methods can be found in the literature for the numerical integration of
Eq. (\ref{gen.eq}), it is important to remark that if the order $r \ge 3$, then necessarily some of the coefficient
$a_j$ and $b_j$ have to be negative \cite{blanes05otn,sheng89spd,suzuki90fdo}. This, while does not constitute a
particular problem when the differential equation is reversible, makes unfeasible their application in parabolic differential
equations of evolutionary type, when the operators $A$ and $B$ are only assumed to generate $C^0$ semi-groups (and not groups): in that case
the flows $\e^{t A}$ and/or $\e^{t B}$ may not be defined for $t < 0$ \cite{castella09smw,hansen09esf,hansen09hos}. Notice that this is the case, in particular, 
if $A$ is the Laplacian operator.

Moreover, even in problems where splitting methods of order $r \ge 3$ can be safely applied,
the presence of negative coefficients usually leads to large truncation errors, so that more stages than strictly necessary
to achieve a given order have to be included in the composition to reduce these errors and improve the overall efficiency \cite{blanes08sac}. 

It is with the aim of circumventing these drawbacks that splitting methods
with complex coefficients (with positive real part) have entered into the literature, mainly in the context of the integration of parabolic differential equations
\cite{castella09smw,hansen09hos,blanes13oho},
but also for ordinary differential equations (ODEs) when structure-preserving (symplecticity, energy conservation, reversibility) is at stake \cite{chambers03siw}. 

Splitting and composition 
methods with complex coefficients, although computationally between 2 and 4 times more costly than their real counterparts when applied to ODEs
involving real vector fields,  possess however
some remarkable properties: their truncation errors with the minimum number of stages are typically very small, and their stability threshold is comparatively
large. Moreover, when the numerical solution is projected at each time step, they lead to approximations that still preserve important
qualitative features (such as symplecticity and time-symmetry) up to an
order much higher than the order of the method itself \cite{blanes10smw,casas21cop,blanes21osc}. 


To better illustrate these points, let us consider a time-symmetric second order method $\mathcal{S}(h)$ (such as one of the compositions (\ref{strang}). Then,
a fourth-order method can be obtained by composition. More specifically, since the coefficients of such a scheme have to satisfy three order conditions,  it makes sense to take three maps,
\[
    \mathcal{S}(\gamma_3 h) \, \mathcal{S}(\gamma_2 h) \, \mathcal{S}(\gamma_1 h).
\]  
In that case, the order conditions read \cite{blanes08sac,hairer06gni}
\begin{equation}  \label{orcon1}
	\sum_{j=1}^3 \gamma_j=1, \qquad   
	\sum_{j=1}^3 \gamma_j^3=0,  \qquad  
	\sum_{j=1}^{2} \left( \gamma_j^3 \left( \sum_{k=j+1}^3 \gamma_k \right) -  \gamma_j \left( \sum_{k=j+1}^3 \gamma_k^3 \right) \right)=0.
\end{equation}
and admit only one real solution, namely
\[
   \gamma_1 =\gamma_3 = \frac{1}{2 - 2^{1/3}}
	, \qquad \gamma_2 = 1- 2 \gamma_1
\]
leading to a time-symmetric composition scheme, usually referred to as Yoshida's method, here denoted as $\mathcal{S}^{[4]}(h)$. 
Notice, however, that there are four more complex solutions. The first pair,
\begin{equation}  \label{palin}
   \gamma_1=\gamma_3 \equiv \gamma = \frac{1}{2 - 2^{1/3} \e^{2 i k \pi/3}}
	, \qquad
	\gamma_2 = 1- 2 \gamma
	, \qquad  k=1,2
\end{equation}
leads again to two time-symmetric methods, denoted as $\Psi_{P,c}^{[4]}(h)$, whereas the second one, denoted as $\Psi_{SC,c}^{[4]}(h)$,
\begin{equation} \label{sym-conj}
	\gamma_1=\bar \gamma_3=  \frac{1}4 \pm i \, \frac14\sqrt{\frac53}
	,  	\qquad \gamma_2 = \frac12  
\end{equation}	
(here the bar indicates the complex conjugate), corresponds to a so--called  \emph{symmetric-conjugate} 
composition method \cite{blanes21osc}: it is symmetric in the real part of the coefficients and
skew-symmetric in the imaginary part. {Here and in the sequel, the first sub-index in a method (either $P$ or $SC$) refers to  its type (either palindromic or
symmetric-conjugate, respectively), whereas the second sub-index (either $r$ or $c$) indicates that the $a_i$ coefficients in the splitting are real or complex,
respectively.}

At order five there are two additional order conditions. One of them, $ \omega_{5,1} = \sum_{j=1}^3 \gamma_j^5$, has been typically used to measure the relative error of methods of the same class. If one defines the error as ${\cal E}=|\omega_{5,1}|$, then one has for the previous methods the following values of
${\cal E}$:
\[
\begin{array}{ccc}
  \mathcal{S}^{[4]}(h) &  \ & {\cal E}=5.29\ldots, \\
  \Psi_{P,c}^{[4]}(h) &  \  & {\cal E}=0.024\ldots, \\
  \Psi_{SC,c}^{[4]}(h) & \ & {\cal E}=0.027\ldots 
\end{array}  
\]
 Notice that the error of methods with complex coefficients is about 200 smaller than in the real case.
 

In the particular case in which $\mathcal{S}(h)$ is given by (\ref{strang}), the previous methods can also be written as
\begin{equation} \label{yoshida4}
 \e^{b_4 h B} \, \e^{a_3 h A} \, \e^{b_3 h B} \, \e^{a_2 h A} \,  \e^{b_2 h B} \, \e^{a_1 h A} \, \e^{b_1 h B},
\end{equation}
with  
\[
b_1 = \frac12\gamma_1, \  \  \  a_1 = \gamma_1, \  \   \ 
b_2 = \frac12(\gamma_1+\gamma_2), \  \  \  a_2 = \gamma_2, \  \  \ 
b_3 = \frac12(\gamma_2+\gamma_3), \  \  \  a_3 = \gamma_3, \  \  \ 
b_4 = \frac12\gamma_3.
\]
As a matter of fact, the simplest symmetric-conjugate composition corresponds to the third order scheme
\begin{equation}  \label{compo3}
  \Psi_{SC,c}^{[3]}(h) =  \mathcal{S}_{\alpha_2 h}^{[2]} \circ    \mathcal{S}_{\alpha_1 h}^{[2]},
\end{equation}
with
\[
   \alpha_1 = \bar{\alpha}_2 \equiv \alpha = \frac{1}{2} + i \frac{\sqrt{3}}{6}.
\]   
For equation (\ref{evol.1}), method (\ref{compo3}) can be written as 
\begin{equation} \label{split31}
   \Psi_{SC,c}^{[3]}(h) = \e^{\bar{b}_1 h B} \, \e^{\bar{a}_1 h A} \, \e^{b_2 h B} \, \e^{a_1 h A} \,  \e^{b_1 h B},
\end{equation}
with $a_1 = \alpha$, $b_1 = \alpha/2$, $b_2 = 1/2$. 

Although (\ref{split31}) is of order 3, if $A$ and $B$ are real, then it renders a scheme of order 4
when it is projected on the real axis after each time step. In addition, it verifies 
$ \Psi_{SC,c}^{[3]}(-h)  \circ \Psi_{SC,c}^{[3]}(h) = I + \mathcal{O}(h^8)$. It is said that the scheme is pseudo-symmetric of order 7, {since
it preserves the time-symmetry property up to terms of order $h^7$} \cite{casas21cop}.


Schemes with complex coefficients have been proposed before for the treatment of quantum problems, 
mainly in the context of imaginary time propagation, with the purpose of computing ground state energies \cite{bandrauk06cis}
and in quantum Monte Carlo simulations \cite{suzuki91gto,goth20hoa}, but also in the decomposition of unitary operators \cite{prosen06hon}. 
In the later case it is shown, both for
unitary $2 \times 2$ matrices and empirically for exponentials of Gaussian random Hermitian matrices, that a splitting method does indeed possess a maximal
time step for which the scheme is numerically stable. We generalize the treatment to differential equations defined in 
$\mathrm{SU}(2)$ for methods possessing a particular
symmetry and eventually examine their behavior when they are applied to the time dependent Schr\"odinger equation.

\section{Splitting methods in $\mathrm{SU}(2)$}

In the study of the evolution of two-level quantum systems one has to deal with the Schr\"odinger equation, which in this context reads ($\hbar = 1$)
\begin{equation} \label{u2.1}
 i  \frac{dU}{dt} = H \, U, \qquad U(0) = I,
\end{equation}
where $U(t)$ is a $2 \times 2$ unitary matrix with unit determinant and the skew-Hermitian Hamiltonian $H$ can be expressed as a linear combination of Pauli
matrices,
\begin{equation} \label{pauli}
 \sigma_1 =
       \left( \begin{array} {cc}
             0    &   1  \\
                 1    &   0
         \end{array}  \right) , \qquad
 \sigma_2 =
       \left( \begin{array} {cr}
             0    &  -i  \\
                 i    &   0
         \end{array}  \right) , \qquad
 \sigma_3 =
       \left( \begin{array} {cr}
             1    &   0  \\
                 0    &  -1
         \end{array}  \right).
\end{equation}
Since our purpose is to analyze splitting methods in this context, we assume that $H$ can be written as
\begin{equation} \label{ham1}
   H = \mathbf{a} \cdot \mathbf{\sigma} + \mathbf{b} \cdot \mathbf{\sigma}
\end{equation}
for given vectors $\mathbf{a}, \mathbf{b} \in \mathbb{R}^3$ and $\mathbf{\sigma} = (\sigma_1, \sigma_2, \sigma_3)$, so that, by comparing with (\ref{evol.1}), 
one has $A \equiv -i \, \mathbf{a} \cdot \mathbf{\sigma}$ and $B \equiv   -i \, \mathbf{b} \cdot \mathbf{\sigma}$.
The exact solution of
Eq. (\ref{u2.1}) after one time step $h$ is 
\[
  U_{\mathrm{ex}}(h) = \e^{-i h H} = \e^{-i h (\mathbf{a} \cdot \mathbf{\sigma} + \mathbf{b} \cdot \mathbf{\sigma})}.
\]
On the other hand, if a splitting method of the form (\ref{td.2b}) of order $r$ with \emph{real} coefficients is applied to solve this very simple problem, it is clear that
the corresponding approximation can be written as 
\[
  U_{\mathrm{app}}(h) = \e^{-i \, h \, \mathbf{d}(h) \cdot \mathbf{\sigma}}, \qquad \mbox{ where } \qquad   \mathbf{d}(h) = \mathbf{a} + \mathbf{b} +
    \mathcal{O}(h^{r}) \in \mathbb{R}^3,
\]    
and thus the method still renders an approximation in $\mathrm{SU}(2)$. The situation is different, however, when the splitting method (\ref{td.2b})
involves \emph{complex} coefficients, since in that case $\mathbf{d}(h) \in \mathbb{C}^3$ and the approximation is no longer unitary. In general, the scheme will be
unstable and the errors will grow exponentially along the integration.

\paragraph{Example.}
At this point it is worth testing the previous third- and fourth-order schemes obtained by composing the Strang splitting and involving complex
coefficients, namely $\Psi_{SC,c}^{[3]}(h)$, $\Psi_{SC,c}^{[4]}(h)$, and $\Psi_{P,c}^{[4]}(h)$. To do that, we consider the following simple Hamiltonian in  
 $\mathrm{SU}(2)$: $H = \sigma_1+\sigma_2$, or alternatively, $\mathbf{a}=(1,0,0)$, $\mathbf{b}=(0,1,0)$ in (\ref{ham1}).

\begin{figure}[!ht] 
\centering
  \includegraphics[width=.49\textwidth]{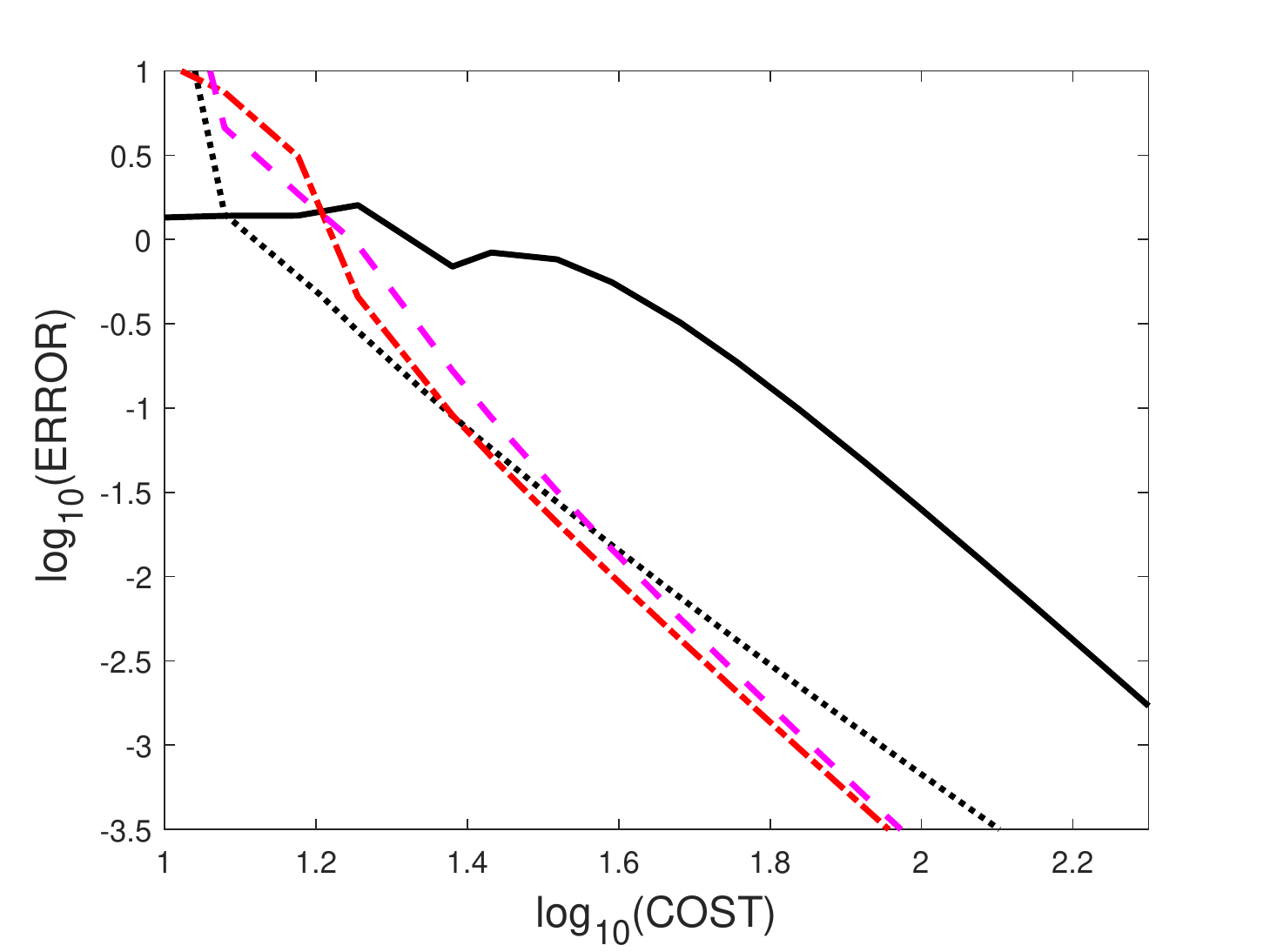}
  \includegraphics[width=.49\textwidth]{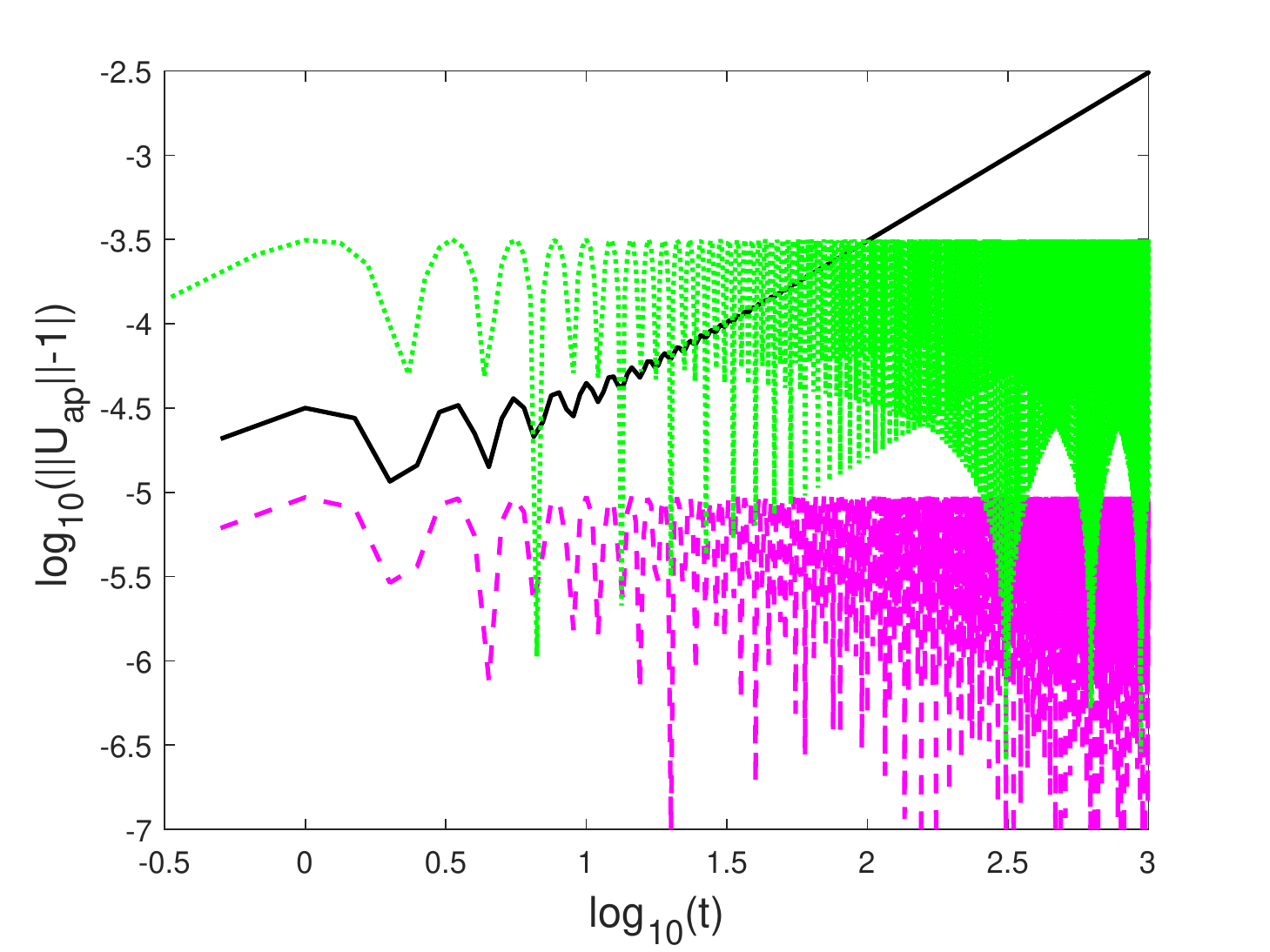}
\caption{\label{Figure0} \small Left: 2-norm error vs. computational cost (number of exponentials) for $\Psi_{SC,c}^{[3]}$ (dotted line), $\mathcal{S}^{[4]}$ (real coefficients, solid line), $\Psi_{P,c}^{[4]}$ (complex coefficients, dash-dotted line) and $\Psi_{SC,c}^{[4]}$ (dashed line) . Right: Error in unitarity for $\Psi_{P,c}^{[4]}$  
(solid line) and the symmetric-conjugate methods  $\Psi_{SC,c}^{[3]}$ (dotted line) and $\Psi_{SC,c}^{[4]}$ (dashed line).
}
\end{figure}

In our experiment, we take as initial condition $U(t_0=0)=I$,  integrate Eq. (\ref{u2.1}) with different values of the time step $h$ and compute
the error of the approximation (in the 2-norm) at the final time $t_f = 10$ as a function of the computational cost (estimated as the number of exponentials 
involved in the whole integration). The results obtained with each method
are displayed in Figure~\ref{Figure0} (left panel). We notice that all schemes involving complex coefficients provide considerably more accurate results than
$\mathcal{S}^{[4]}(h)$ (black solid line), the fourth-order methods being also more efficient than $\Psi_{SC,c}^{[3]}(h)$ for high accuracy. 

In order to check how each scheme with complex coefficients behaves with respect to unitarity, we take as a final time $t_f=1000$, and adjust $h$ (and therefore 
the number of iterations $N$) so that 
they require the same computational cost. Specifically, $N=6000$ ($h=1/6$) for scheme  $\Psi_{SC,c}^{[3]}(h)$ and $N=4000$ ($h=1/4$) for all methods of
order 4.
Finally we compute $|\|U_{\mathrm{app}}(nh)\|-1|, \ n=1,2,\ldots, N=(t_f-t_0)/h$, where $U_{\mathrm{app}}(nh)$ denotes the approximate solution after $n$ steps. 
The outcome is depicted in Figure~\ref{Figure0} (right panel). Notice how the error in unitarity grows for $\Psi_{P,c}^{[4]}(h)$, whereas it is bounded, even for
large intervals, for the symmetric-conjugate methods  $\Psi_{SC,c}^{[3]}(h)$ and $\Psi_{SC,c}^{[4]}(h)$. Among them, the later clearly provides more accurate results.

\begin{figure}[!ht] 
\centering
  \includegraphics[width=.6\textwidth]{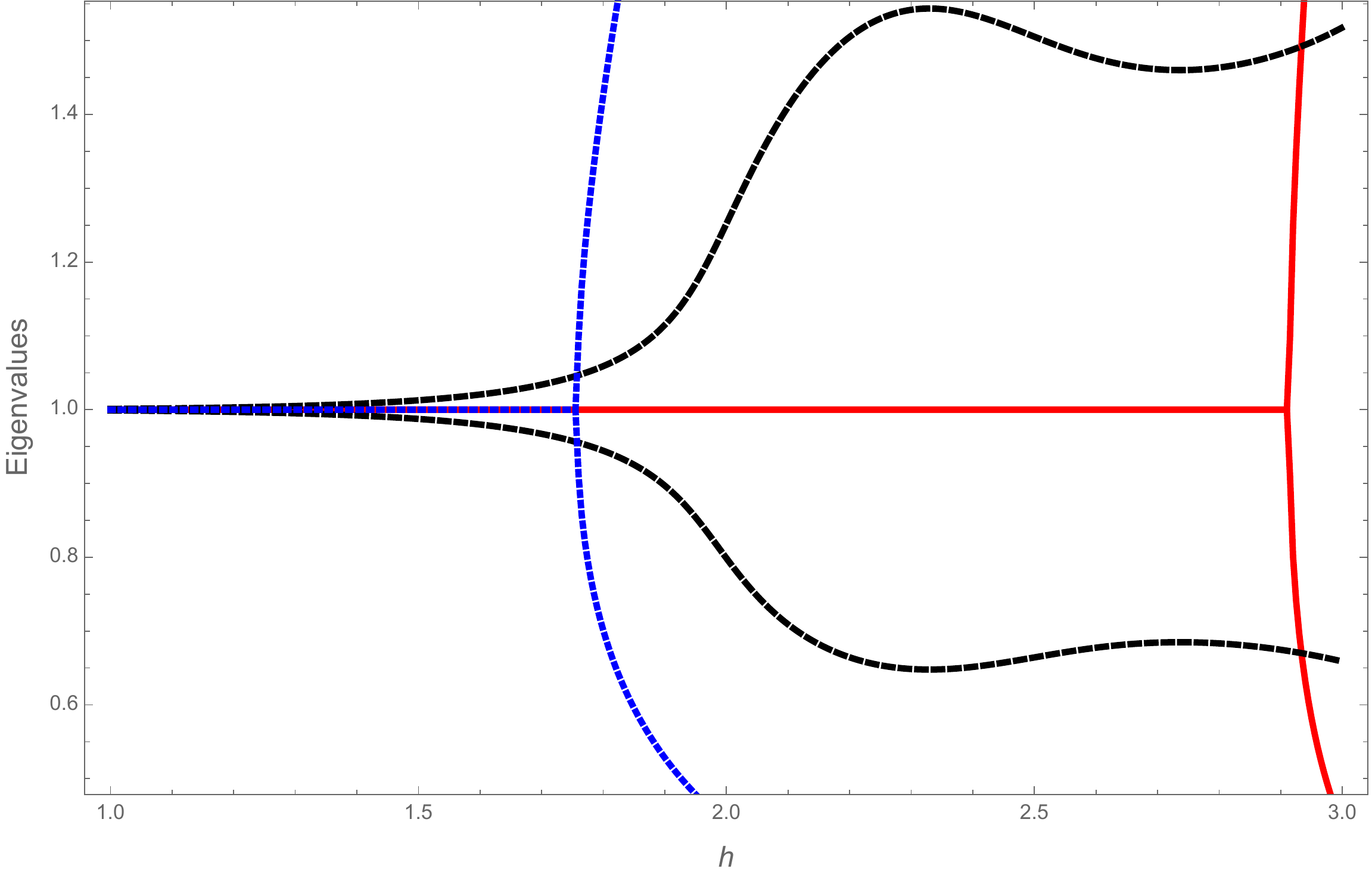}
\caption{\label{Figure0a} \small Absolute value of the eigenvalues of the approximate solution matrix obtained with $\Psi_{P,c}^{[4]}$ with complex coefficients ($k=1$, black dashed line), $\Psi_{SC,c}^{[3]}$ (blue dotted line) and $\Psi_{SC,c}^{[4]}$ (red, solid line).
}
\end{figure}

This marked difference of both types of integrators can also be illustrated by computing the eigenvalues $\lambda_1$, $\lambda_2$
of the approximate solution after one step size, i.e.,
of the corresponding matrix $U_{\mathrm{app}}(h)$. In the exact case, of course, both evolve on the unit circle in the complex plane, whereas here 
one has still $\lambda_1 \lambda_2 = 1$, since the determinant is one.
In Figure~\ref{Figure0a} we depict $|\lambda_j|$, $j=1,2$, as a function of $h$ for the palindromic scheme $\Psi_{P,c}^{[4]}$ with $k=1$ (black, dashed lines)
and the symmetric-conjugate splittings $\Psi_{SC,c}^{[3]}$ (blue, dotted line) and $\Psi_{SC,c}^{[4]}$ (red, solid line) in the range $1 \le h \le 3$. It is worth remarking that for the symmetric-conjugate methods both $|\lambda_j|$ are
exactly 1 for $0 \le h \le h^*$, with  $h^*=1.7570473$ for $\Psi_{SC,c}^{[3]}$ and $h^*= 2.9139468357$ for $\Psi_{SC,c}^{[4]}$. In other words, they behave
as unitary maps when $h \le h^*$. On the other hand, it can be checked that $|\lambda_1| > 1$ for any $h>0$ for $\Psi_{P,c}^{[4]}$.  $\Box$

\

The previous example illustrates in fact a general pattern  exhibited by symmetric-conjugate methods for this problem, as we next prove.


{
\begin{proposition}  \label{propo.1}
Suppose a symmetric-conjugate splitting method of the form (\ref{td.2b}), with $a_{s+1-j} = \bar{a}_j$, $b_{s+2-j} = \bar{b}_j$, is applied to the numerical
integration of the Schr\"odinger equation (\ref{u2.1}) with the Hamiltonian given by (\ref{ham1}). In that case, the following statements hold:
\begin{itemize}
  \item[(a)] The eigenvalues of the matrix approximating
the solution after on time step $h$ lie on the unit circle in the complex plane for sufficiently small $h$.
  \item[(b)] The symmetric-conjugate splitting method is itself conjugate to a unitary method for sufficiently small $h$.
 \end{itemize} 
\end{proposition}
}
\begin{Proof}
When a splitting method of the form (\ref{td.2b}) is applied to solve Eq. (\ref{u2.1}),
the corresponding approximation after one step can be written as $U_{\mathrm{app}}(h) = \exp(V(h))$, where $V(h)$ is a linear combination of $A$, $B$ and all their
nested commutators,
\begin{equation} \label{modif1}
  V(h) = h(w_{1,1} A+ w_{1,2} B) + h^2 w_{2,1} [A,B] + h^3 (w_{3,1} [A,[A,B]] + w_{3,2} [B,[A,B]]) + \mathcal{O}(h^4),
\end{equation}  
and $w_{n,k}$ are polynomials in the coefficients $a_j$, $b_j$. Method (\ref{td.2b}) is of order $r$ if $w_{1,1} = w_{1,2} = 1$ and
the polynomials $w_{n,k}$ vanish for $1 < n < r$. In our case, since
\[
[-i \, \mathbf{a} \cdot \mathbf{\sigma}, -i \, \mathbf{b} \cdot \mathbf{\sigma}] = - i \, 2 \, (\mathbf{a} \times \mathbf{b}) \cdot \mathbf{\sigma},
\]
it is clear that the vector fields associated with all commutators in (\ref{modif1}) containing an even number of operators are perpendicular to the plane generated by the vectors $\mathbf{a}$ and $\mathbf{b}$, whereas those containing an odd number 
of operators $A$ and $B$ are in such a plane. If in addition the method is symmetric-conjugate, then a straightforward computation shows that
$V^{\dag}(h) = V(-h)$, and so all polynomials $w_{2j+1,k}$ are real whereas  all polynomials $w_{2j,k}$ are pure imaginary. Therefore,
$V(h)$ can be written as 
\begin{equation} \label{modif2}
  V(h) = - i \, h \tilde{H}(h), \qquad \mbox{ with } \qquad \tilde{H}(h) = \mathbf{d}(h) \cdot \mathbf{\sigma} + i \, \mathbf{c}(h) \cdot \mathbf{\sigma},
\end{equation}
for two vectors $\mathbf{c}, \mathbf{d} \in \mathbb{R}^3$ verifying $\mathbf{c} \cdot \mathbf{d} = 0$ and
\begin{equation} \label{modif3}
   \mathbf{d}(h) = \mathbf{a} + \mathbf{b} + \mathcal{O}(h^{r+1}), \qquad\qquad \mathbf{c}(h) =  \mathcal{O}(h^{r}).
\end{equation}
This special structure of $\tilde{H}(h)$ allows one to obtain statements (a) and (b) above. First, if we write 
\[ 
  \e^{- i \, h \tilde{H}(h)} = \e^{h \, \mathbf{u} \cdot \mathbf{\sigma}}, \qquad \mbox{ with } \qquad \mathbf{u} = \mathbf{c} - i \, \mathbf{d},
\]
then 
\[
     U_{\mathrm{app}}(h) = \e^{h \, \mathbf{u} \cdot \mathbf{\sigma}} = \cosh( h u) I + \frac{\sinh(h u)}{u} \, \mathbf{u} \cdot \mathbf{\sigma},
\]
with $u = \sqrt{ \mathbf{u} \cdot \mathbf{u}} = (\|\mathbf{c}\|^2 - \|\mathbf{d}\|^2)^{1/2}$. Of course, if  $\|\mathbf{c}\| <  \|\mathbf{d}\|$, then 
$\cosh(u) = \cos \alpha$, $\sinh(u) = i \sin \alpha$, with $\alpha = (\|\mathbf{d}\|^2 - \|\mathbf{c}\|^2)^{1/2}$, and the eigenvalues of 
$U_{\mathrm{app}}(h)$ are $\lambda_{1,2} = \exp(\pm i h \alpha(h))$. But, in virtue of (\ref{modif3}),  this always holds for sufficiently small values
of $h$. 

Statement (b) can demonstrated as follows. Let us  introduce the unitary vector
\[
 \mathbf{C} = \frac{\mathbf{d} \times \mathbf{c}}{\|\mathbf{d} \times \mathbf{c}\|}.
\]
A trivial computation shows that
\[
   \mathbf{d} \times \mathbf{C} = - \frac{\|\mathbf{d}\|}{\|\mathbf{c}\|} \, \mathbf{c}, \qquad\qquad
   \mathbf{c} \times \mathbf{C} =   \frac{\|\mathbf{c}\|}{\|\mathbf{d}\|} \, \mathbf{d},
\]
and furthermore, for a given parameter $s \in \mathbb{R}$,
\[
  \e^{s \, \mathbf{C} \cdot \mathbf{\sigma}} \, \e^{h \mathbf{u} \cdot \mathbf{\sigma}} \, \e^{-s \, \mathbf{C} \cdot \mathbf{\sigma}} = \exp \left(
   \e^{s \, \mathbf{C} \cdot \mathbf{\sigma}} \, (h \, \mathbf{u} \cdot \mathbf{\sigma}) \, \e^{-s \, \mathbf{C} \cdot \mathbf{\sigma}} \right).
\]
From the definition of $\mathbf{C}$ and the properties of the Pauli matrices \cite{galindo90qme}, one has
\begin{eqnarray*}
   \e^{s \, \mathbf{C} \cdot \mathbf{\sigma}} \, (\mathbf{c} \cdot \mathbf{\sigma}) \, \e^{-s \, \mathbf{C} \cdot \mathbf{\sigma}}  & = &
   \cosh(2 s) (\mathbf{c} \cdot \mathbf{\sigma}) - i \, \sinh(2 s)  \frac{\|\mathbf{c}\|}{\|\mathbf{d}\|} \, \mathbf{d} \cdot \mathbf{\sigma} \\
   \e^{s \, \mathbf{C} \cdot \mathbf{\sigma}} \, (\mathbf{d} \cdot \mathbf{\sigma}) \, \e^{-s \, \mathbf{C} \cdot \mathbf{\sigma}}  & = &
   \cosh(2 s) (\mathbf{d} \cdot \mathbf{\sigma}) + i \, \sinh(2 s)  \frac{\|\mathbf{d}\|}{\|\mathbf{c}\|} \, \mathbf{c} \cdot \mathbf{\sigma}, 
\end{eqnarray*}
and thus
\[
   \e^{s \, \mathbf{C} \cdot \mathbf{\sigma}} \,  (\mathbf{u} \cdot \mathbf{\sigma}) \, \e^{-s \, \mathbf{C} \cdot \mathbf{\sigma}} =
   \left( \sinh(2 s) \frac{\|\mathbf{d}\|}{\|\mathbf{c}\|} + \cosh(2 s) \right) \mathbf{c} \cdot \mathbf{\sigma} - i 
   \left( \sinh(2 s) \frac{\|\mathbf{c}\|}{\|\mathbf{d}\|} + \cosh(2 s) \right) \mathbf{d} \cdot \mathbf{\sigma}.
\]
If we now take $s$ such that
\begin{equation} \label{tanh}
 \sinh(2 s) \frac{\|\mathbf{d}\|}{\|\mathbf{c}\|} + \cosh(2 s) = 0, \qquad \mbox{ i.e., } \qquad \tanh(2 s) = -\frac{\|\mathbf{c}\|}{\|\mathbf{d}\|},
\end{equation}
then, clearly
\[
   \e^{s \, \mathbf{C} \cdot \mathbf{\sigma}} \,  (-i h \tilde{H}) \, \e^{-s \, \mathbf{C} \cdot \mathbf{\sigma}} =
    - i \, h \, \mathbf{D}(h)  \cdot \mathbf{\sigma},
\]
with
\[
      \mathbf{D}(h) = \sinh(2 s) \left(  \frac{\|\mathbf{c}\|}{\|\mathbf{d}\|} -  \frac{\|\mathbf{d}\|}{\|\mathbf{c}\|} \right)  \mathbf{d} \in \mathbb{R}^3.
\]
In consequence,
\begin{equation} \label{conju1}
        \e^{s \, \mathbf{C} \cdot \mathbf{\sigma}} \,  \e^{- i h \tilde{H}} \, \e^{-s \, \mathbf{C} \cdot \mathbf{\sigma}} =
        \e^{ - i h \mathbf{D}(h)  \cdot \mathbf{\sigma}}.
\end{equation}
In other words, if $s$ is such that Eq. (\ref{tanh}) holds, then the map obtained by applying a symmetric-conjugate splitting method 
is conjugate to a unitary matrix. Notice that if $\|\mathbf{c}\| < \|\mathbf{d}\|$ this is always possible, in agreement with statement (a)
for the eigenvalues of the approximate solution matrix.   
\end{Proof}   

{
Proposition \ref{propo.1} thus provides a rigorous justification of the results shown in Figures \ref{Figure0} and \ref{Figure0a}: since a
symmetric-conjugate splitting method is ultimately conjugate to a unitary map in the sense of eq. (\ref{conju1}) for sufficiently small
values of $h$, then the error in the
unitarity of the numerical solution is bounded, whereas the eigenvalues remain on the unit circle in the complex plane.
}

{
Methods of the form (\ref{conju1}) are called \emph{processed} or \emph{corrected} in the literature 
(see, e.g. \cite{blanes16aci,blanes08sac,hairer06gni,mclachlan02sm}). In that context, method
$ \e^{- i h \tilde{H}}$ is called the kernel, and $ \e^{s \, \mathbf{C} \cdot \mathbf{\sigma}}$ the processor. For integrators of this class, only
the error terms in the kernel that cannot be removed by a processor are relevant in the long run. In the case of unitary problems in
$\mathrm{SU}(2)$ we have shown that any symmetric-conjugate splitting method is indeed the kernel of a processed unitary scheme.
}

\section{Application to the time-dependent Schr\"odinger equation}

In view of the previous results in $\mathrm{SU}(2)$, it is natural to examine the situation when splitting methods with complex coefficients, and in particular
symmetric-conjugate schemes, are applied in a more general setting. To this end, we next consider 
the numerical integration of the general time dependent Schr\"odinger equation
\begin{equation}   \label{Schr0}
  i \frac{\partial}{\partial t} \psi (x,t)  = -\frac{1}{2\mu}
  \Delta \psi (x,t) + V(x) \psi (x,t),
\end{equation}
where now
$\psi:\mathbb{R}^d\times\mathbb{R} \longrightarrow \mathbb{C}$ is the wave function representing the state of the system and the
initial state is $\psi(x,0)=\psi_0(x)$. We take again $\hbar = 1$ and a reduced mass $\mu$.
Equation (\ref{Schr0}) can be written as
\begin{equation}   \label{Schr1}
  i \frac{\partial}{\partial t} \psi   = (\hat{T}(P) + \hat{V}(X)) \psi,
\end{equation}
with $\hat{T}(P) = P^2/(2 \mu)$, and the operators $X$ and $P$ are defined by their
actions on $\psi(x,t)$ as
\begin{equation} \label{xp}
  X \psi (x,t) = x\, \psi (x,t), \qquad\quad
  P \ \psi (x,t) = -i \, \nabla \psi (x,t).
\end{equation}
The usual procedure for applying splitting methods in this setting consists first in 
discretizing the space variables $x$, so as to get  a system of 
ordinary differential equations (ODEs) which is subsequently integrated in time by the splitting scheme. If, for simplicity, we consider the one-dimensional problem,
$d=1$, and suppose that it is defined
in  $x \in [x_0, x_N]$, first this interval is partitioned
 into $N$ parts of length $\Delta x = (x_{N}-x_{0})/N$ and the vector $u = (u_{0}, \ldots,
u_{N-1})^T \in \mathbb{C}^N$ is formed, with $u_{n} =
\psi(x_{n},t)$ and $x_{n} = x_{0} + n \Delta x$,
$n=0,1,\ldots,N-1$. The partial differential equation
(\ref{Schr0}) is then replaced by the $N$-dimensional linear ODE
\begin{equation}   \label{lin1}
  i \frac{d }{dt} u(t) = H \, u(t),  \qquad
    u(0) = u_{0} \in \mathbb{C}^N,
\end{equation}
where now $H$ represents the (real symmetric) $N \times N$ matrix associated with the Hamiltonian. 

When a Fourier spectral collocation method is used, then the matrix $H$ in (\ref{lin1}) is 
\begin{equation} \label{ham1}
  H = T + V, 
\end{equation}  
where $V$ is a diagonal matrix associated with
the potential $\hat{V}$ and $T$ is a (full) differentiation matrix related with the kinetic energy $\hat{T}$. Their action on the wave function vector $u$
is trivial: on the one hand, $(V u)_{n} = V(x_{n}) u_{n}$ and thus the product
$V  u$ requires to compute $N$ complex multiplications. On the other hand, 
$T u = \mathcal{F}^{-1} D_T
\mathcal{F} u$, where $\mathcal{F}$ and $\mathcal{F}^{-1}$
are the forward and backward discrete Fourier transform,
and $D_T$ is again diagonal. The transformation $\mathcal{F}$ from the
discrete coordinate representation to the discrete momentum
representation (and back) is done via the fast Fourier transform
(FFT) algorithm, requiring $\mathcal{O}(N \log N)$ operations.

Notice that, since 
\[
\left(\e^{\tau V} u \right)_i =
\e^{\tau V(x_i)}u_i, \qquad\quad 
\e^{\tau T} u = \mathcal{F}^{-1} \e^{\tau D_T} \mathcal{F} \, u,
\] 
splitting methods constitute a valid alternative to approximate the solution $u(t) = \e^{ \tau (T+V)} u_0$ for a time step $\Delta t$, with $\tau = -i \Delta t$. 
Thus, with the Lie--Trotter scheme (\ref{lie-trotter}) one has
\[
  \e^{\tau (T+V)} = \e^{\tau T} \, \e^{\tau V} + \mathcal{O}(\tau^2),
\]
whereas the 2nd-order Strang splitting (\ref{strang}) constructs the numerical approximation $u_{n+1}$ at time $t_{n+1} = t_n + \Delta t$ by
\[
  u_{n+1} =   \e^{\tau/2 V} \, \e^{\tau T} \, \e^{\tau/2 V} \, u_n \equiv \mathcal{S}(\tau) u_n.
\]
The resulting scheme is called split-step Fourier method in the chemical literature, and has some remarkable properties. In particular, it is both unitary and symplectic
\cite{lubich08fqt}, as well as time-reversible. In addition, for suitable regularity assumptions on the potential and on the norm of the commutators 
$[\hat{T}, \hat{V}] = \hat{T} \hat{V} - \hat{V} \hat{T}$ and
$[\hat{T}, [\hat{T}, \hat{V}]]$, the error at $t_n$ is bounded by 
\[
  \| u_n - u(t) \| \le C \, \Delta t^2 \, t \, \max_{0 \le s \le t} \|u(s)\|_2.
\]  
Higher order methods can be obtained of course  by considering compositions (\ref{td.2b}), which in this setting read
\begin{equation}   \label{td.2}
 \Psi^{[r]}(\tau) = \e^{b_{s+1} \tau V} \e^{a_s \tau T}  \e^{b_s \tau V} \cdots  \e^{a_1 \tau T}  \, \e^{b_1 \tau V},
\end{equation}
and in fact, a large collection of practical schemes of different orders exist for carrying out the numerical integration
(see e.g. \cite{blanes08sac,hairer06gni,mclachlan02sm} and references therein). In addition, 
from (\ref{xp}), 
it is clear that $ [\hat{X},\hat{P}] \ \psi (x,t) = i \, \psi (x,t),$ and so
\begin{equation} \label{rkn}
   [\hat{V},[\hat{V},[\hat{V},\hat{T}]]] \ \psi (x,t) = 0.
\end{equation}
This property leads to a reduction in the number of order conditions necessary to achieve a given order $r$ and allows one to construct highly
efficient schemes.



\section{Splitting methods with complex coefficients}
\label{sec.complex}

When exploring the applicability of splitting methods with complex coefficients to the general time-dependent Schr\"odinger equation, 
several aspects must be addressed. First, since the computational cost of method (\ref{td.2}) is dominated by the number of FFTs per step, 
the presence of complex $a_j$, $b_j$ does not contribute significantly to increase this cost. In addition, it has been shown in other problems that
splitting methods with complex coefficients
involving the minimum number of flows to achieve a given order already provide good efficiency, in contrast with their real counterparts.
On the other hand, however, since $\sum_j a_j = 1$ for a consistent method, if $a_j \in \mathbb{C}$, then 
imaginary parts positive \emph{and} negative enter into the game, with the result that severe instabilities may arise in practice due to the unboundedness of the
Laplace operator. With respect to the potential, since in regions where it takes large values
the wave function typically is close to zero, we can introduce an artificial cut-off bound in the computation if necessary, so that complex $b_j$ can in principle be used, at least for a sufficiently small $\Delta t$.
It makes sense therefore to construct and examine in detail methods with real $a_j$ and complex $b_j$ coefficients.

In the following, and for simplicity, we restrict ourselves to splitting methods (\ref{td.2}) of order $r \le 4$, which we denote by their sequence of coefficients
as
\[
   (b_{s+1},a_s, b_s, \ldots, a_2, b_2, a_1, b_1).
\]   
The order conditions are then
\begin{eqnarray} \label{orcon}
 \mbox{Order 1:} & \quad  &  \sum_{i=1}^s a_i = 1, \qquad\qquad\qquad\qquad   \sum_{i=1}^s b_i = 1, \nonumber \\
  \mbox{Order 2:} & \quad  &  \sum_{i=1}^s b_i \left( \sum_{j=1}^i a_j \right) = \frac{1}{2}, \nonumber \\
  \mbox{Order 3:} & \quad  &  \sum_{i=1}^s b_i \left( \sum_{j=1}^i a_j \right)^2 = \frac{1}{3}, \qquad\quad  
     \sum_{i=1}^s a_i \left( \sum_{j=i}^s b_j \right)^2 = \frac{1}{3},  \\
  \mbox{Order 4:} & \quad  &  \sum_{i=1}^s b_i \left( \sum_{j=1}^i a_j \right)^3 = \frac{1}{4}, \qquad\quad  
     \sum_{i=1}^s a_i \left( \sum_{j=i}^s b_j \right)^3 = \frac{1}{4}, \nonumber \\
    & \quad  &  \sum_{i=1}^s a_i^2  \left( \sum_{j=i}^s b_j \right)^2 + 2 \, \sum_{i=2}^s a_i \left( \sum_{j=1}^{i-1} a_j \right)
      \left( \sum_{k=i}^s b_k \right)^2 = \frac{1}{6}. \nonumber
\end{eqnarray}      
In typical applications of splitting methods with real coefficients, 
only palindromic sequences of coefficients, i.e., methods (\ref{td.2}) with $b_{s+2-j} = b_j$, $a_{s+1-j}=a_j$ for all $j$ 
are used. In that case, all the conditions at even order are automatically satisfied and the resulting schemes are time-symmetric, 
$\Psi^{[r]}(\tau) \, \Psi^{[r]}(-\tau) = I$, and of even order. Here, however, since we are dealing with complex coefficients, we also
analyze the case $r=3$ for completeness.

\paragraph{Order 3.} The first five order conditions in (\ref{orcon}) admit solutions with all $a_j$ real and positive and $b_j \in \mathbb{C}$ with
positive real part if one considers a composition of the form
\begin{equation} \label{split32}
  \Psi_{SC,r}^{[3]}(\tau) = (\bar{b}_1, a_1, \bar{b}_2, a_2, b_2, a_1, b_1)
\end{equation}
involving 6 parameters. Then one gets a 1-parametric family of solutions (+c.c.) with the required properties. Among them, we choose
\[
 a_1= \frac3{10}, \qquad a_2=\frac2{5}, \qquad b_1=\frac{13}{126}-i\frac{\sqrt{59/2}}{63}, \qquad b_2=\frac{25}{63}+i\frac{5\sqrt{59/2}}{126}. 
\]
Composition (\ref{split31}) constitutes of course another symmetric-conjugate method of order 3, denoted here by
\begin{equation} \label{split33}
   \Psi_{SC,c}^{[3]}(\tau) = (\bar{b}_1, \bar{a}_1, b_2, a_1, b_1)
\end{equation}   
and involving less maps, although in this case $a_1 \in \mathbb{C}$.

\paragraph{Order 4.}  

The simplest approach to construct a palindromic scheme with $a_j \in \mathbb{R}$ and $b_j \in \mathbb{C}$  consists in  taking all the $a_j$
equal. In that case, with $s=4$, one has enough parameters to solve the required four order conditions (at odd orders).
Only two solutions (complex conjugate to each other) are obtained, as shown in
\cite{castella09smw}, thus resulting in the scheme
\begin{equation} \label{split42}
  \Psi_{P,r}^{[4]}(\tau) = (b_1, a_1, b_2, a_2, b_3, a_2, b_2, a_1, b_1)
\end{equation}
with
\[
 a_1=a_2=\frac14, \qquad b_1=\frac1{10}-i\frac{1}{30}, \qquad b_2=\frac4{15}+i\frac{2}{15}, \qquad b_3=\frac4{15}-i\frac{1}{5}.
\]
Although more efficient schemes can be obtained if one allows for different $a_j$'s \cite{blanes13oho}, since we are interested here mainly in the qualitative
behavior of the different methods, we limit ourselves to  (\ref{split42}) as representative of palindromic splitting methods with real $a_j$'s and complex $b_j$'s,
whereas we can take scheme (\ref{yoshida4})
\begin{equation} \label{split42b}
  \Psi_{P,c}^{[4]}(\tau) = (b_1, a_1, b_2, a_2, b_2, a_1, b_1), 
\end{equation}
as representative of palindromic methods with both $a_j \in \mathbb{C}$ and $b_j \in \mathbb{C}$.

Symmetric-conjugate splitting methods with real $a_j$'s require at least $s=5$ stages, in which case one has a free parameter. If we fix this as $a_1 = 1/8$,
we get the scheme
\begin{equation} \label{split43}
  \Psi_{SC,r}^{[4]}(\tau) = (\bar{b}_1, a_1, \bar{b}_2, a_2, \bar{b}_3, a_3, b_3, a_2, b_2, a_1, b_1)
\end{equation}
with
\[
\begin{aligned}
 & a_2 = 0.23670501659941197298,   \\
 & a_3 = 0.27658996680117605403,  \\
 & b_1 = 0.03881396214419327198 -   0.045572109263923104872 \, i, \\
 & b_2 = 0.19047619047619047619 +  0.115462072300408741306 \, i,\\
 & b_3 = 0.27070984737961625182 -    0.148322245509626403888 \, i
\end{aligned}
\] 
It is worth noticing that one can obtain symmetric-conjugate methods from  palindromic schemes and vice versa. Thus, in particular, by
composing the palindromic scheme (\ref{split42}) with its complex conjugate we can form a symmetric-conjugate splitting method with 8 stages
and $a_j \in \mathbb{R}$, $b_j \in \mathbb{C}$:
\begin{equation} \label{sc4}
  \Xi_{SC,r}^{[4]}(\tau) = \Psi_{P,r}^{[4]}(\tau/2)  \, \overline{\Psi}_{P,r}^{[4]}(\tau/2), 
\end{equation}
whereas doing the same with the 3rd-order symmetric-conjugate method (\ref{split32}) results in the 4th-order palindromic scheme with 6 stages and
$a_j \in \mathbb{R}$, $b_j \in \mathbb{C}$:
\begin{equation} \label{p4}
  \Xi_{P,r}^{[4]}(\tau) = \Psi_{SC,r}^{[3]}(\tau/2)  \, \overline{\Psi}_{SC,r}^{[3]}(\tau/2). 
\end{equation}
This is possible because  the adjoint  of $(\Psi_{SC,r}^{[3]}(\tau))^*$ verifies 
\[ 
  (\Psi_{SC,r}^{[3]}(\tau))^* =  \overline{\Psi}_{SC,r}^{[3]}(\tau).
\]      
In our numerical experiments we shall also use for comparison one of the best 4th-order splitting methods with real coefficients designed specifically for systems
verifying (\ref{rkn}). It reads
\begin{equation} \label{rkn4}
 \Psi_{RKN}^{[4]}(\tau) = (b_1, a_1, b_2, a_2, b_3, a_3, b_4, a_3, b_3, a_2, b_2, a_1, b_1) 
\end{equation}
and the coefficients can be found in \cite{blanes02psp}. The scheme has three additional parameters that are used to minimize error terms at higher orders,
and provides by construction unitary approximations.

\section{Numerical experiments}

We next report on some numerical tests we have carried out with the splitting methods presented in section \ref{sec.complex} applied
to the Scr\"odinger equation
in one dimension. Since many different schemes are tested and compared, it is convenient to classify them into the following categories:
\begin{itemize}
 \item symmetric-conjugate methods with $a_j \in \mathbb{R}$, $b_j \in \mathbb{C}$
 \begin{itemize}
  \item Order 3: $\Psi_{SC,r}^{[3]}$, Eq. (\ref{split32});
  \item Order 4:  $\Psi_{SC,r}^{[4]}$, Eq. (\ref{split43});
 \end{itemize}
 \item  symmetric-conjugate  with $a_j \in \mathbb{C}$, $b_j \in \mathbb{C}$: method $\Psi_{SC,c}^{[3]}$, Eq. (\ref{split33}), order 3;
 \item palindromic with $a_j \in \mathbb{R}$, $b_j \in \mathbb{C}$: method $\Psi_{P,r}^{[4]}$, Eq. (\ref{split42}), order 4;
 \item palindromic with  $a_j \in \mathbb{C}$, $b_j \in \mathbb{C}$: method $\Psi_{P,c}^{[4]}$, Eq. (\ref{split42b}), order 4;
\end{itemize}
For completeness, we also consider the following schemes of order 4 with $a_j \in \mathbb{R}$, $b_j \in \mathbb{C}$:
\begin{itemize}
 \item symmetric-conjugate obtained from a palindromic method: $\Xi_{SC,r}^{[4]}$, Eq. (\ref{sc4});
 \item palindromic obtained from a symmetric-conjugate method: $\Xi_{P,r}^{[4]}$, Eq. (\ref{p4}).
\end{itemize}

\paragraph{Quartic potential.}
As the first example we take the quartic oscillator
\begin{equation} \label{quartic}
  V(x) = -\frac{1}{2} x^2 + \frac{1}{20} x^4
\end{equation}
and the initial condition $\psi_0(x) = \sigma \, \e^{-x^2/2}$, where $\sigma$ is a normalization constant.
As usual, and since the exact solution decays rapidly, we truncate the infinite spatial domain to the periodic domain
$[-L, L]$, provided $L$ is sufficiently large and use Fourier spectral methods. We take $L = 8$ and set up a uniform grid on the interval with $N= 128 $ subdivisions.
Finally, we apply the different schemes to integrate in time the resulting equation (\ref{lin1}) in the interval $t \in [0, t_f]$, with $t_f = 8000$. As in the case
of the example in $\mathrm{SU}(2)$, we check the behavior of each scheme with respect to unitarity by computing $|\|u_{\mathrm{app}}(t)\| - 1|$ 
along the integration,
where $u_{\mathrm{app}}(t)$ denotes the numerical approximation obtained by each method.

In addition, we also compute the expected value of the energy, $u_{\mathrm{app}}^*(t) \cdot H  u_{\mathrm{app}}(t)$ and measure the error as the
difference with respect to the exact value:
\begin{equation} \label{eq.5.1b}
   \mbox{energy error:} \quad |u_{\mathrm{app}}^*(t) \cdot (H  u_{\mathrm{app}}(t))  - u_0^* \cdot (H  u_0)|.
\end{equation}   
In each case, the time step is adjusted so that the number of FFTs (and their inverses) are the same for all methods (specifically,
1572864), so that the computational cost of all schemes is similar.

\begin{figure}[!ht] 
\centering
  \includegraphics[width=.49\textwidth]{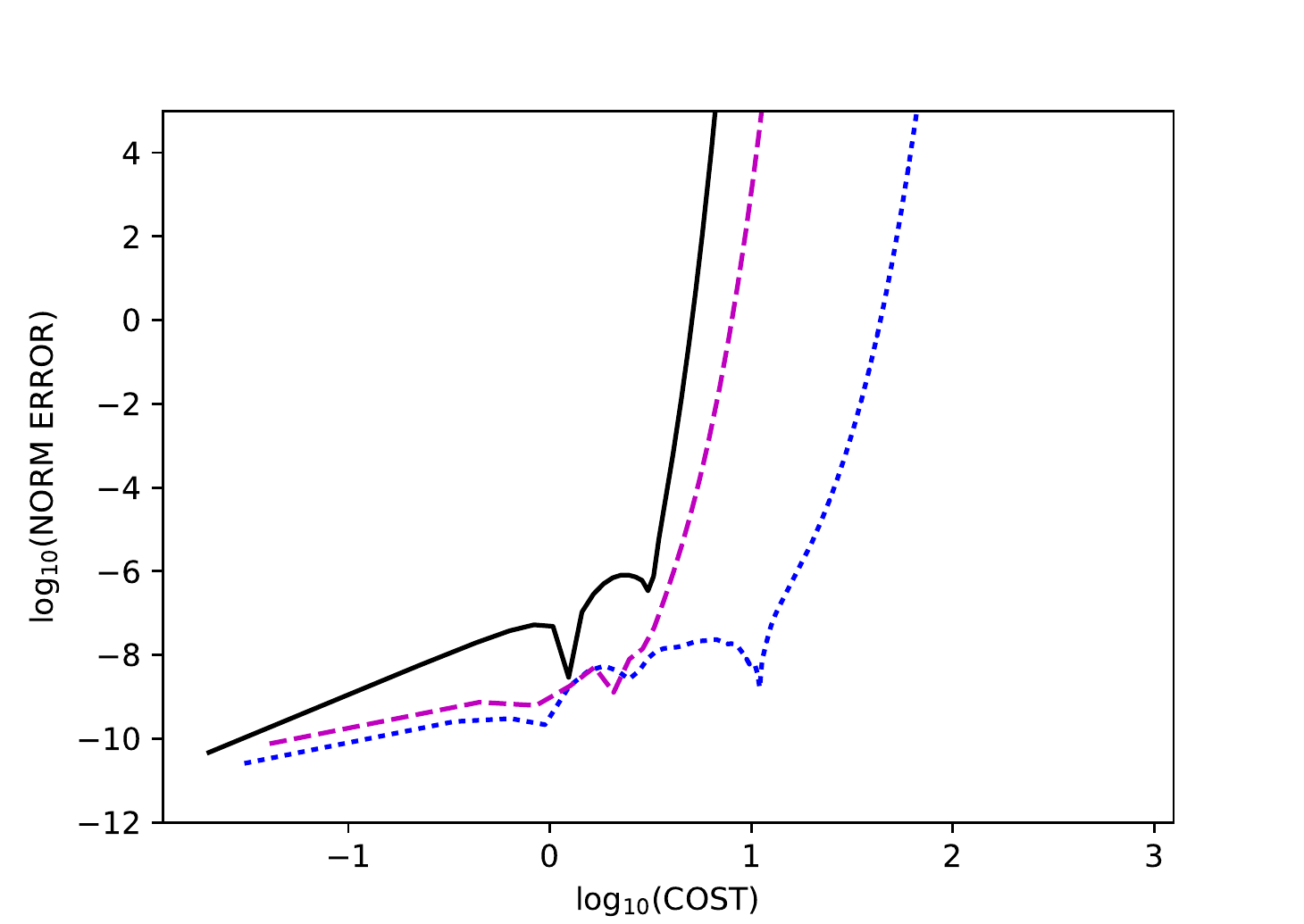}
  \includegraphics[width=.49\textwidth]{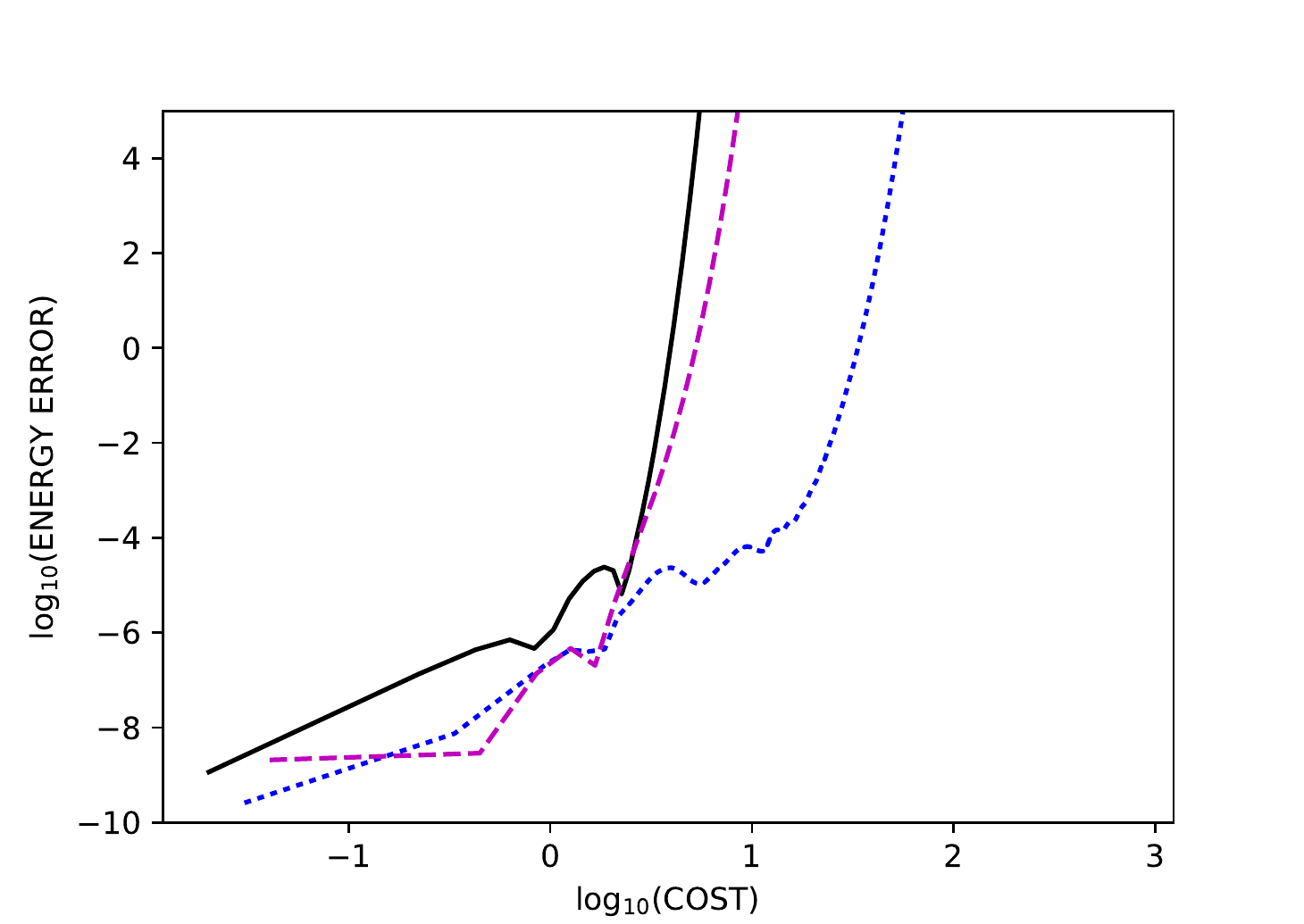}
\caption{\label{figu3} \small Error in norm of the approximate solution (left) and error in energy (\ref{eq.5.1b}) (right) for the quartic potential (\ref{quartic})
obtained by the palindromic schemes $\Psi_{P,r}^{[4]}$ (magenta, dashed line),  $\Psi_{P,c}^{[4]}$ (blue dotted line) and
the symmetric-conjugate method $\Psi_{SC,c}^{[3]}$ (black solid line) along the integration interval. The step size is chosen so that all methods have the same computational cost.
}
\end{figure}

 Figure \ref{figu3} shows the corresponding results obtained by palindromic schemes with the coefficients $a_j$ real, $\Psi_{P,r}^{[4]}$, and $a_j$
 complex, $\Psi_{P,c}^{[4]}$, together with the symmetric-conjugate method $\Psi_{SC,c}^{[3]}$ with $a_j \in \mathbb{C}$. We notice that the
 qualitative behavior of all of them is similar: after some point, depending on the particular step size, the unitarity is lost and the error in energy grows rapidly.

 {
 We notice here the same type of behavior observed in the case of the group $\mathrm{SU}(2)$: palindromic schemes with both real and complex 
 coefficients $a_j$ are unable to preserve unitarity. On the other hand, symmetric-conjugate schemes with $a_j \in \mathbb{C}$ lead also to unbounded
 errors, according with the comments formulated at the beginning of section \ref{sec.complex}.
 }

\begin{figure}[!ht] 
\centering
  \includegraphics[width=.49\textwidth]{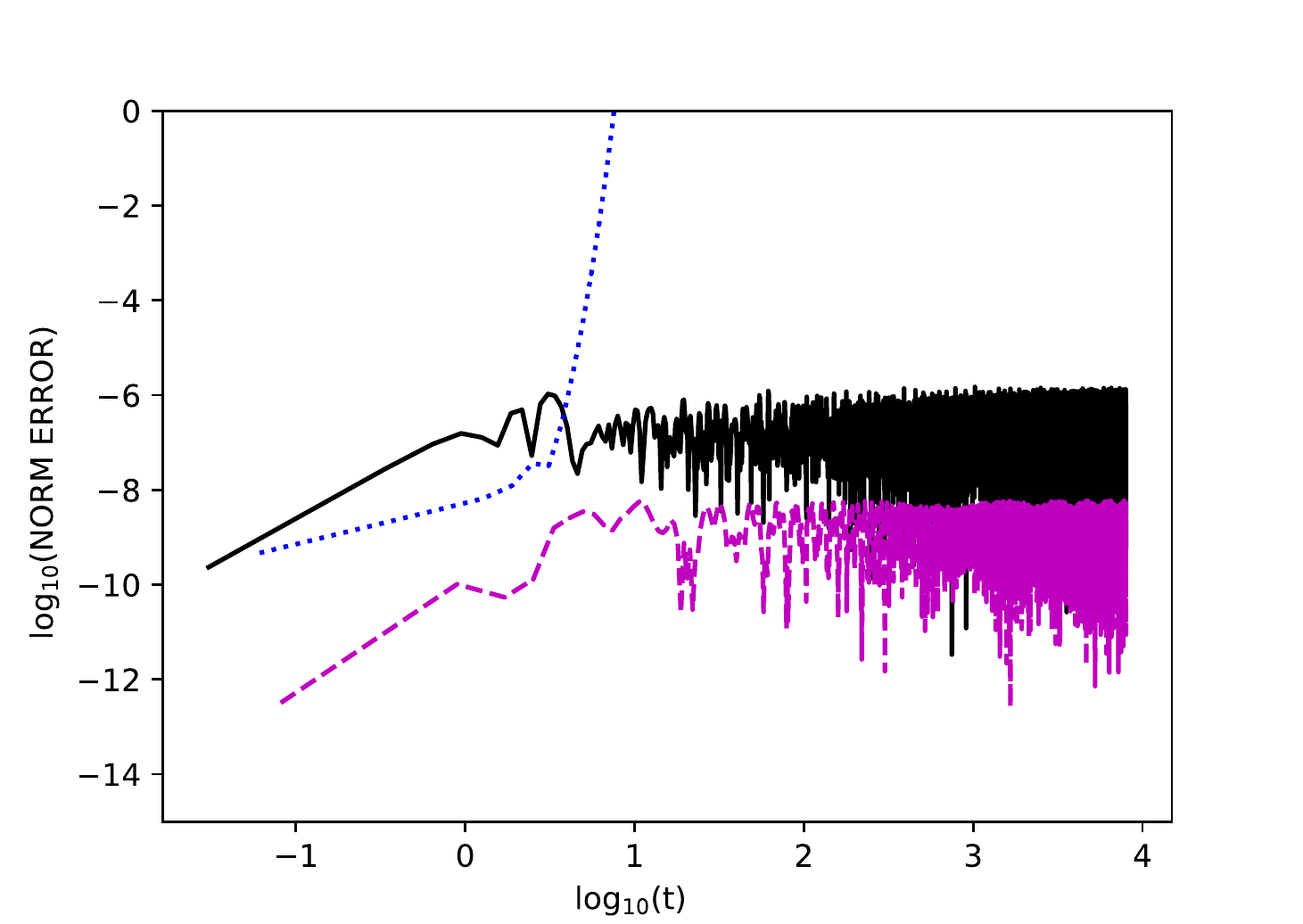}
  \includegraphics[width=.49\textwidth]{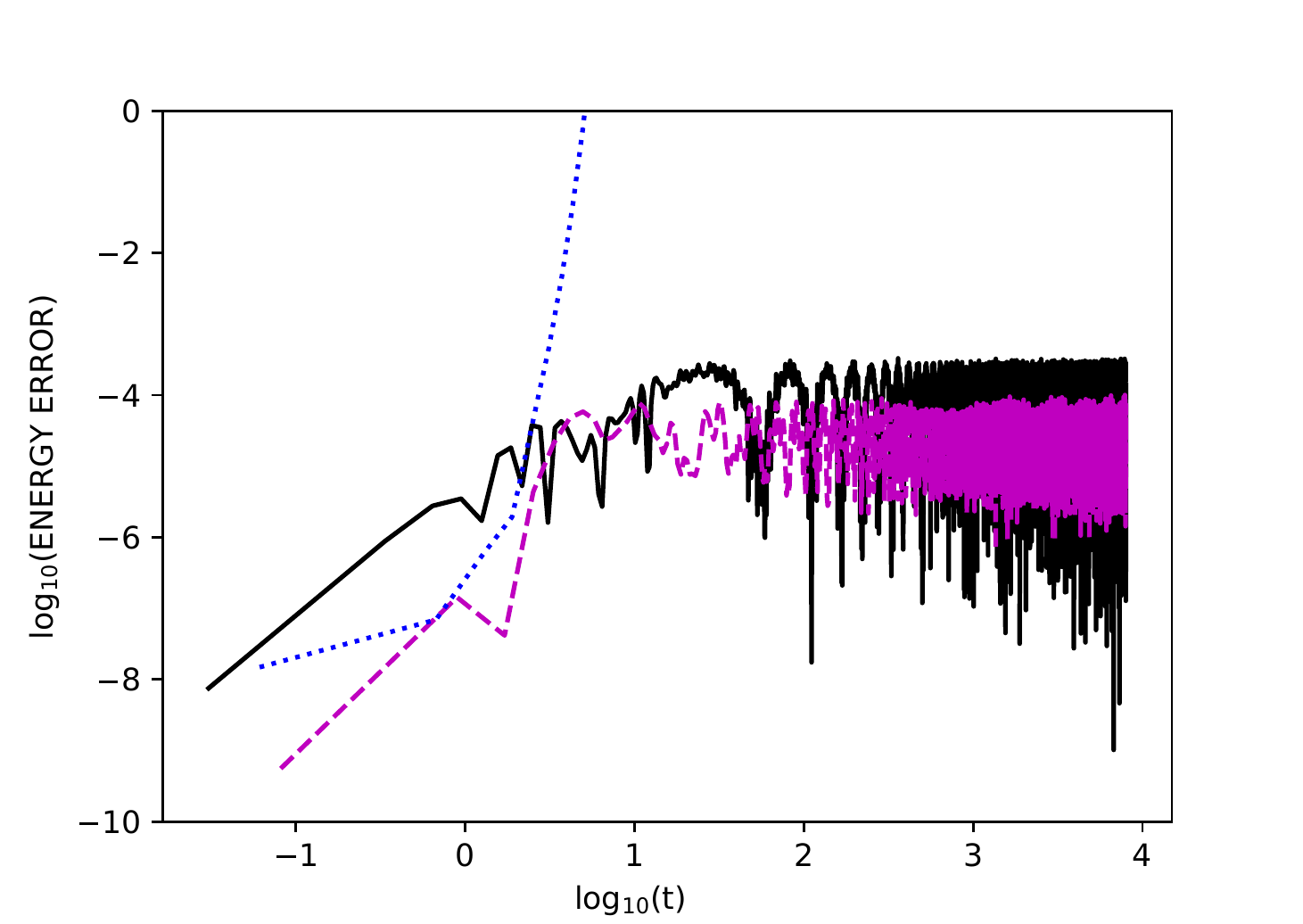}
\caption{\label{figu4} \small Error in norm of the approximate solution (left) and error in energy (\ref{eq.5.1b}) (right) for the quartic potential (\ref{quartic})
obtained by the palindromic scheme $\Xi_{P,r}^{[4]}$ (blue dotted line),  and the symmetric-conjugate
 schemes $\Psi_{SC,r}^{[3]}$ (black solid line) and $\Xi_{SC,r}^{[4]}$ (magenta dashed line) along the integration interval. The step size is chosen so that all methods have the same computational cost.
}
\end{figure}

 We collect in Figure \ref{figu4} the corresponding results achieved by the palindromic method $\Xi_{P,r}^{[4]}$ (blue dotted line), and the symmetric-conjugate
 schemes $\Psi_{SC,r}^{[3]}$ (black solid line) and $\Xi_{SC,r}^{[4]}$ (magenta dashed line), all of them with real parameters $a_j$. It is worth noticing that
 both the norm of the solution and the expected value of the energy are preserved for very long times by symmetric-conjugate methods with $a_j \in \mathbb{R}$,
 and this happens even if the method is obtained by composing a palindromic scheme (with a poor behavior) with its complex conjugate. By contrast,
 a symmetric-conjugate method  looses its good preservation properties when composed to form a palindromic scheme, even if all $a_j$ are real.
 
 {
 We have carried out the same experiment, but with the roles of $T$ and $V$ interchanged. In other words, the complex coefficients $b_j$ are now multiplying
 the discretized Laplacian. In that case, the errors obtained by all the previous schemes grow unbounded. This indicates that, at least for this example,
 one needs both symmetric-conjugate schemes and real coefficients multiplying the Laplacian to get bounded errors in the preservation of unitarity and energy.
 }
 
\begin{figure}[!ht] 
\centering
  \includegraphics[width=.49\textwidth]{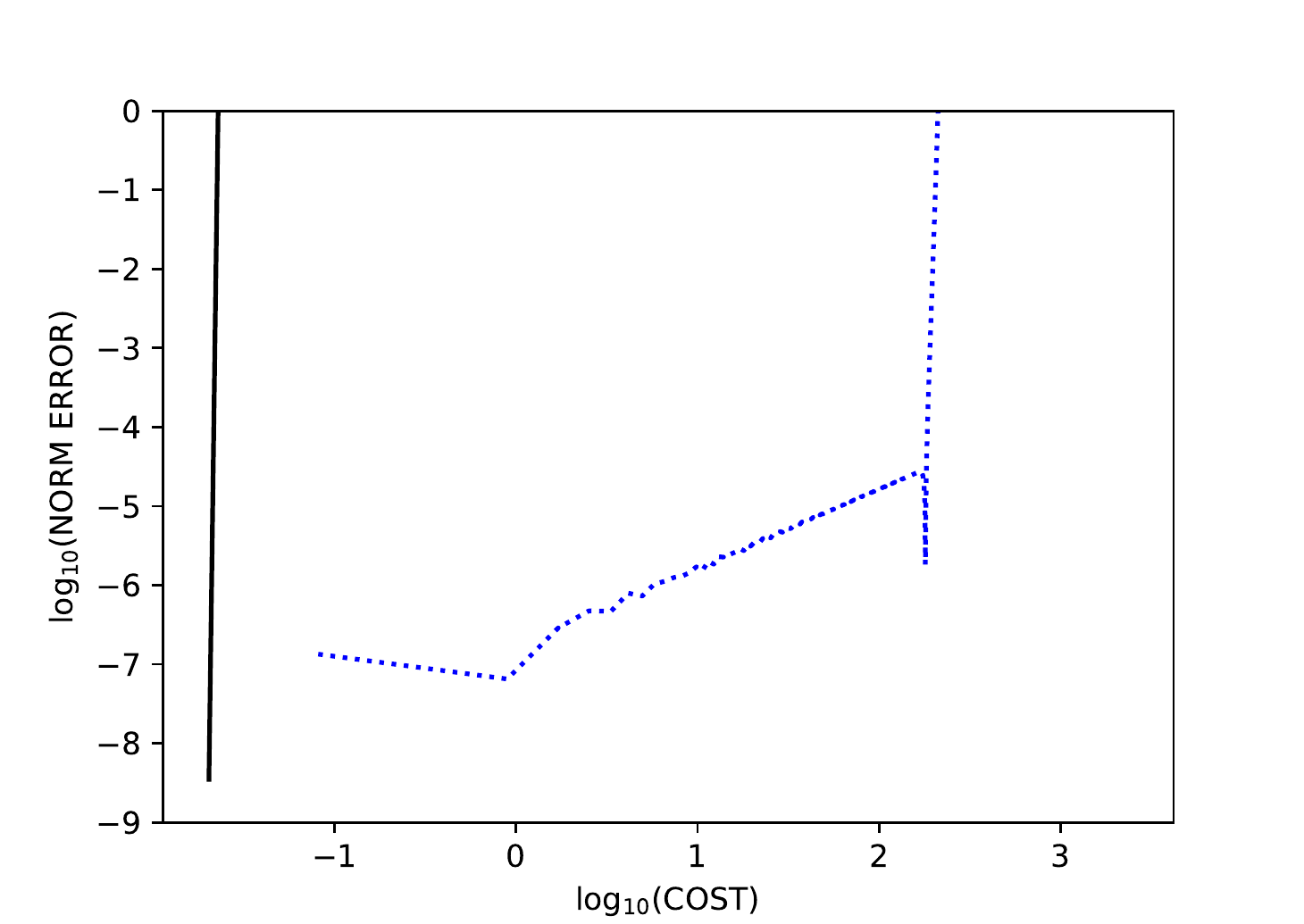}
  \includegraphics[width=.49\textwidth]{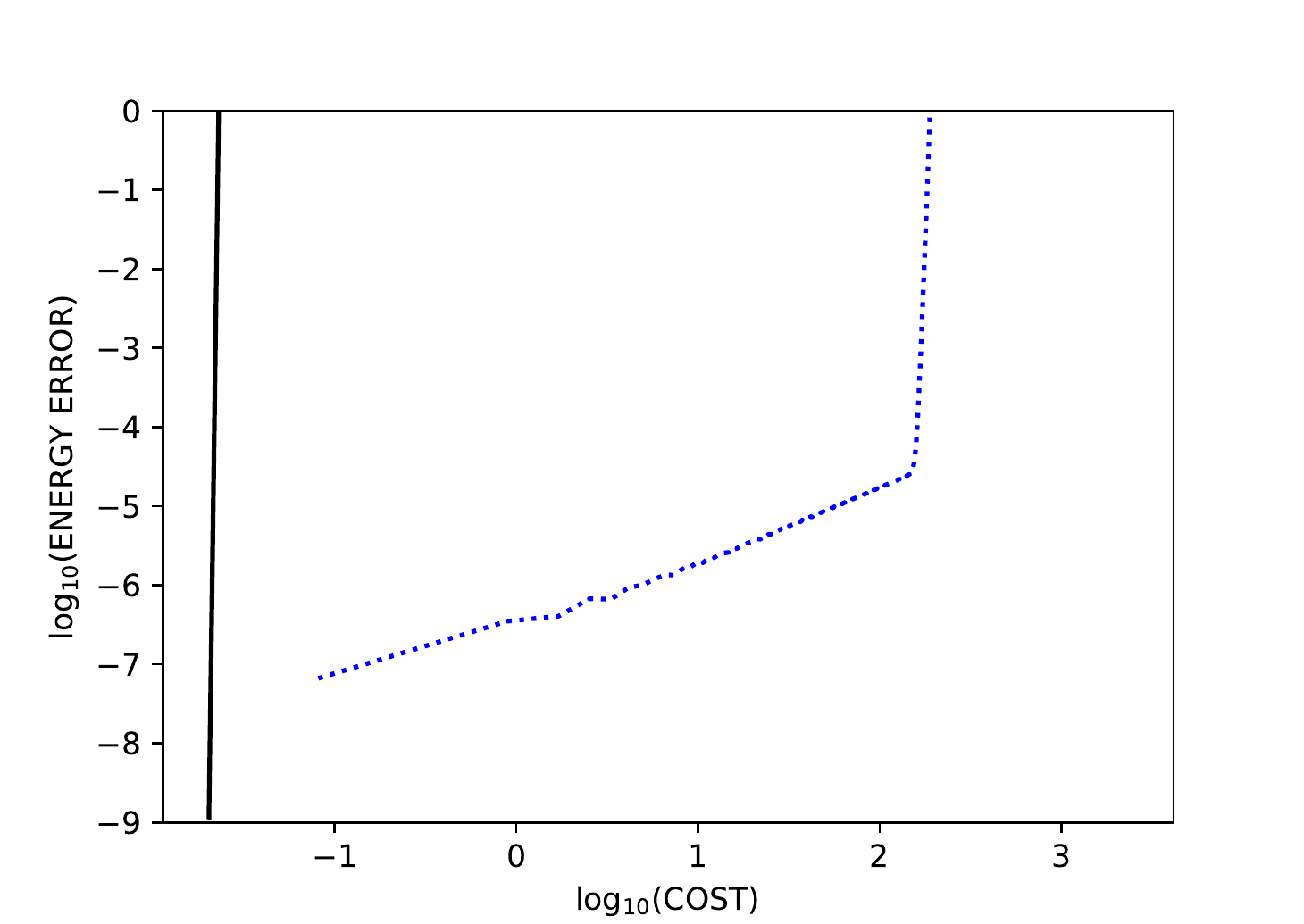}
\caption{\label{figu5} \small Error in norm of the approximate solution (left) and error in energy (\ref{eq.5.1b}) (right) for the P\"oschl--Teller potential (\ref{pt1})
obtained by the palindromic scheme   $\Psi_{P,c}^{[4]}$ (blue dotted line) and
the symmetric-conjugate method $\Psi_{SC,c}^{[3]}$ (black solid line) along the integration interval. The result achieved by $\Psi_{P,r}^{[4]}$ is out of the scale.
}
\end{figure}

\begin{figure}[!ht] 
\centering
  \includegraphics[width=.49\textwidth]{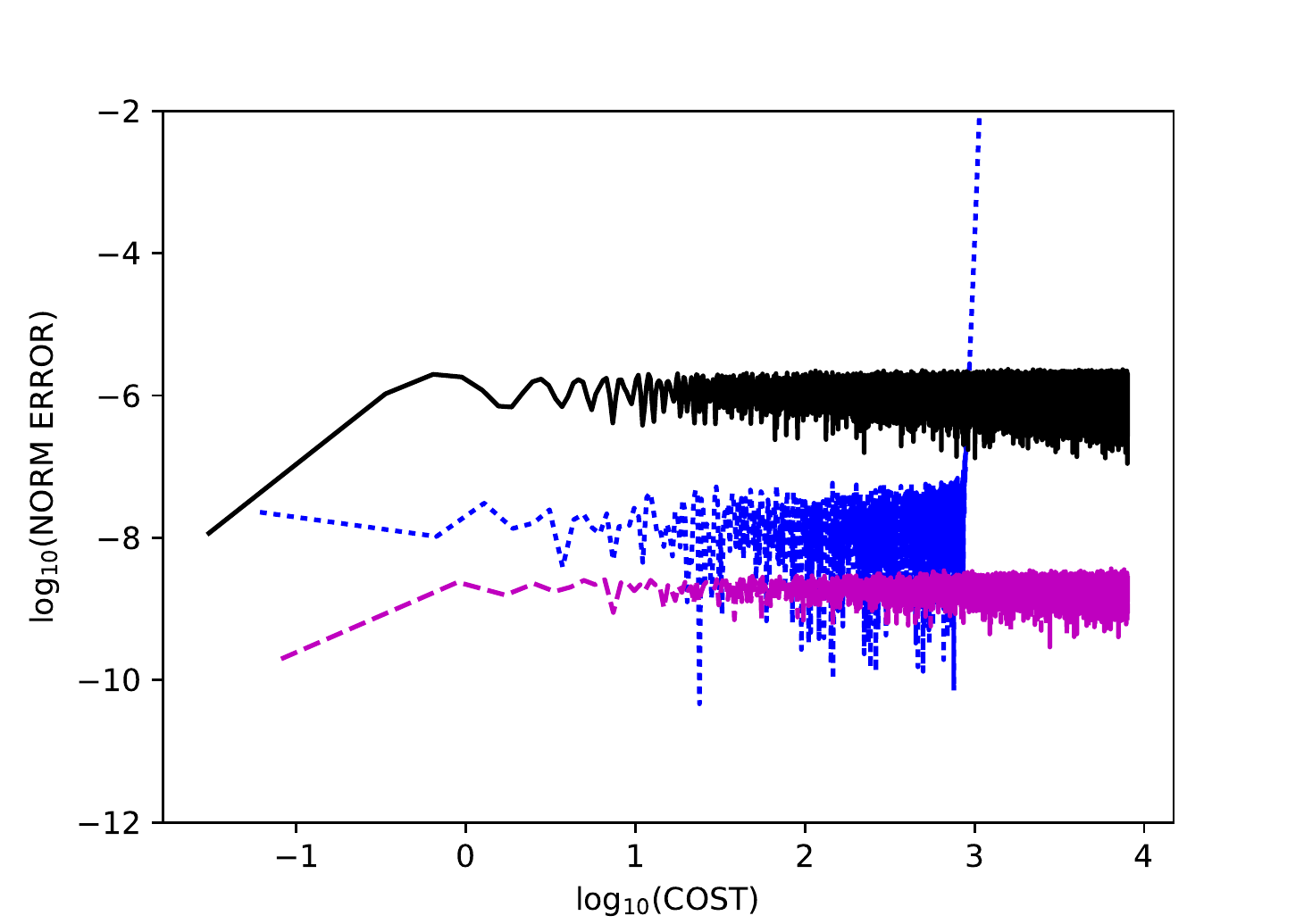}
  \includegraphics[width=.49\textwidth]{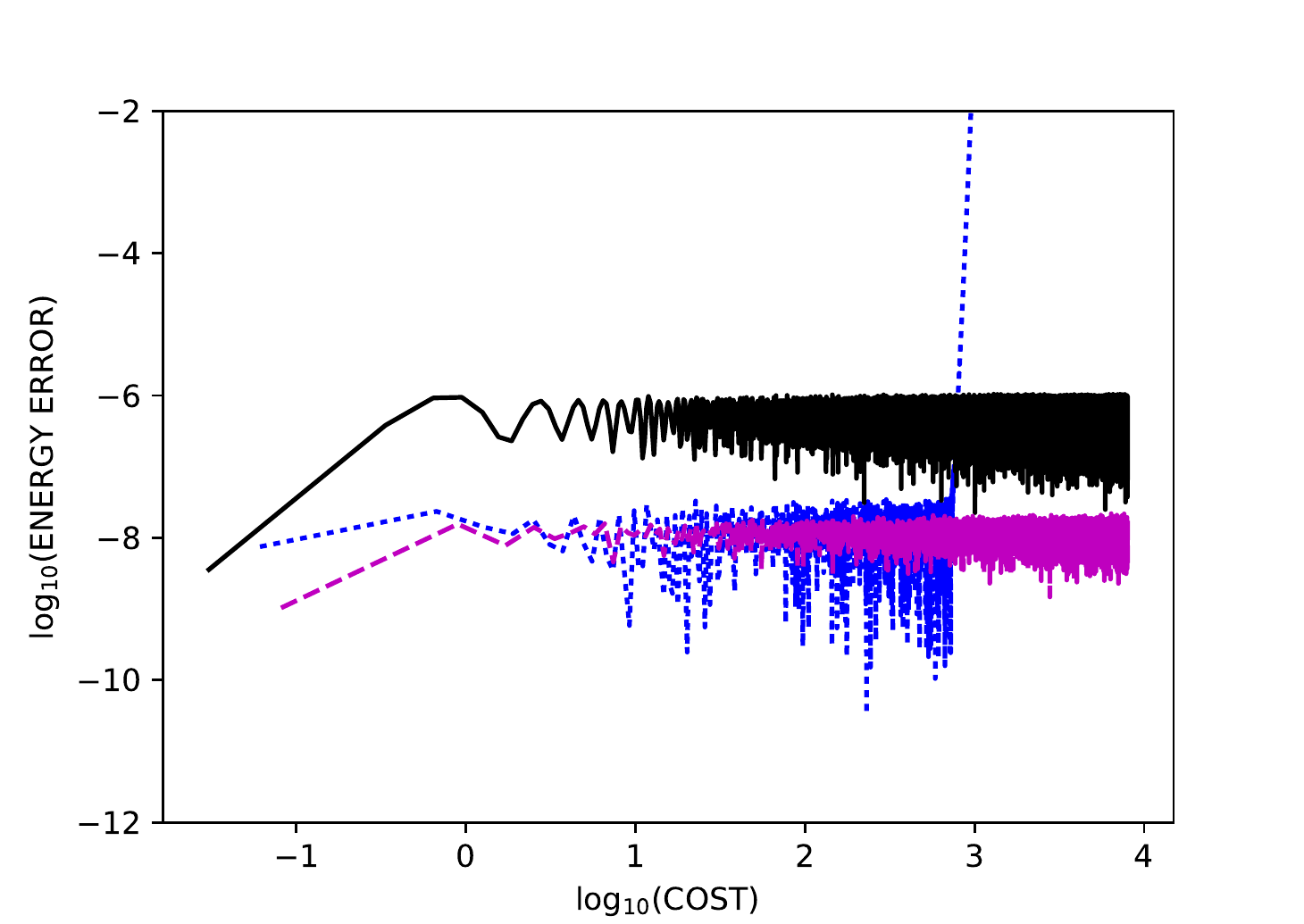}
\caption{\label{figu6} \small Error in norm of the approximate solution (left) and error in energy(\ref{eq.5.1b})  (right) for the P\"oschl--Teller potential (\ref{pt1})
obtained by the palindromic scheme $\Xi_{P,r}^{[4]}$ (blue dotted line),  and the symmetric-conjugate
 schemes $\Psi_{SC,r}^{[3]}$ (black solid line) and $\Xi_{SC,r}^{[4]}$ (magenta dashed line) along the integration interval. 
}
\end{figure}

\paragraph{P\"oschl--Teller potential.} The next set of simulations is carried out  with the well known one-dimensional P\"oschl--Teller potential,
\begin{equation} \label{pt1}
  V(x) = -\frac{\lambda (\lambda+1)}{2} \mbox{sech}^2(x),
\end{equation}
with $\lambda(\lambda+1) = 10$. It has been used in polyatomic molecular simulations and admits an analytic treatment \cite{flugge71pqm}. We take again
as initial condition $\psi_0(x) = \sigma \, \e^{-x^2/2}$, with $\sigma$ a normalizing constant, then apply Fourier spectral methods on the interval 
$x \in [-8,8]$ and integrate until the final time $t_f=8000$ with the previous numerical splitting methods. 
{ For this potential we take $N=512$ subdivisions of the space interval to better visualize the behavior of the different methods. 
Figure  \ref{figu5} is the analogous of Fig. \ref{figu3}), and only displays the results obtained by $\Psi_{SC,c}^{[3]}$ (black solid line) and
$\Psi_{P,c}^{[4]}$ (blue dotted line), since the output corresponding to $\Psi_{P,r}^{[4]}$ is out of the scale  (the errors are  greater than $10^{87}$).  
On the other hand, Figure  \ref{figu6} shows the same pattern as Figure \ref{figu4}: only symmetric-conjugate schemes with $a_j \in \mathbb{R}$ provide
bounded errors in the norm and in the energy of the solution.}


\begin{figure}[!ht] 
\centering
  \includegraphics[width=.49\textwidth]{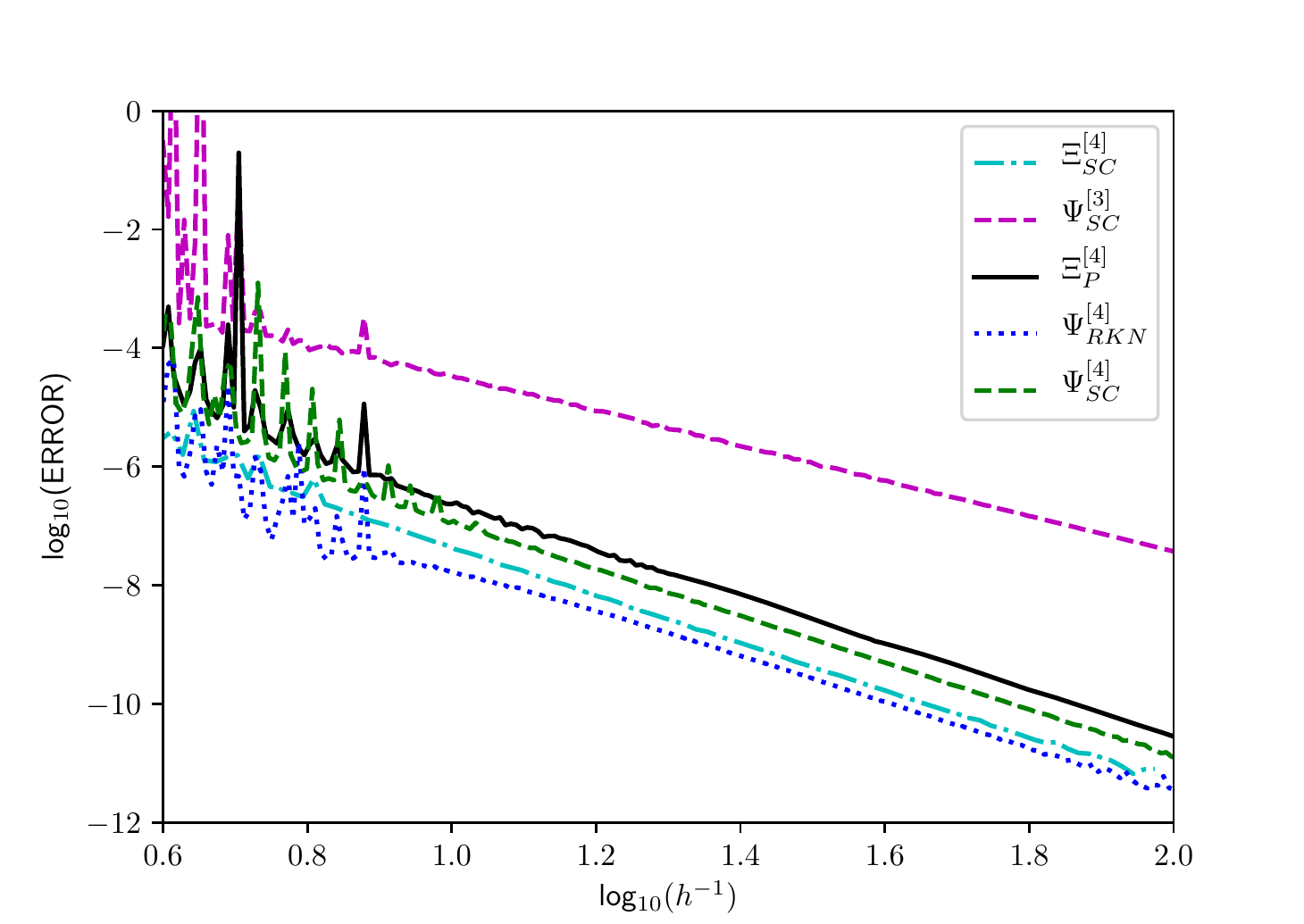}
  \includegraphics[width=.49\textwidth]{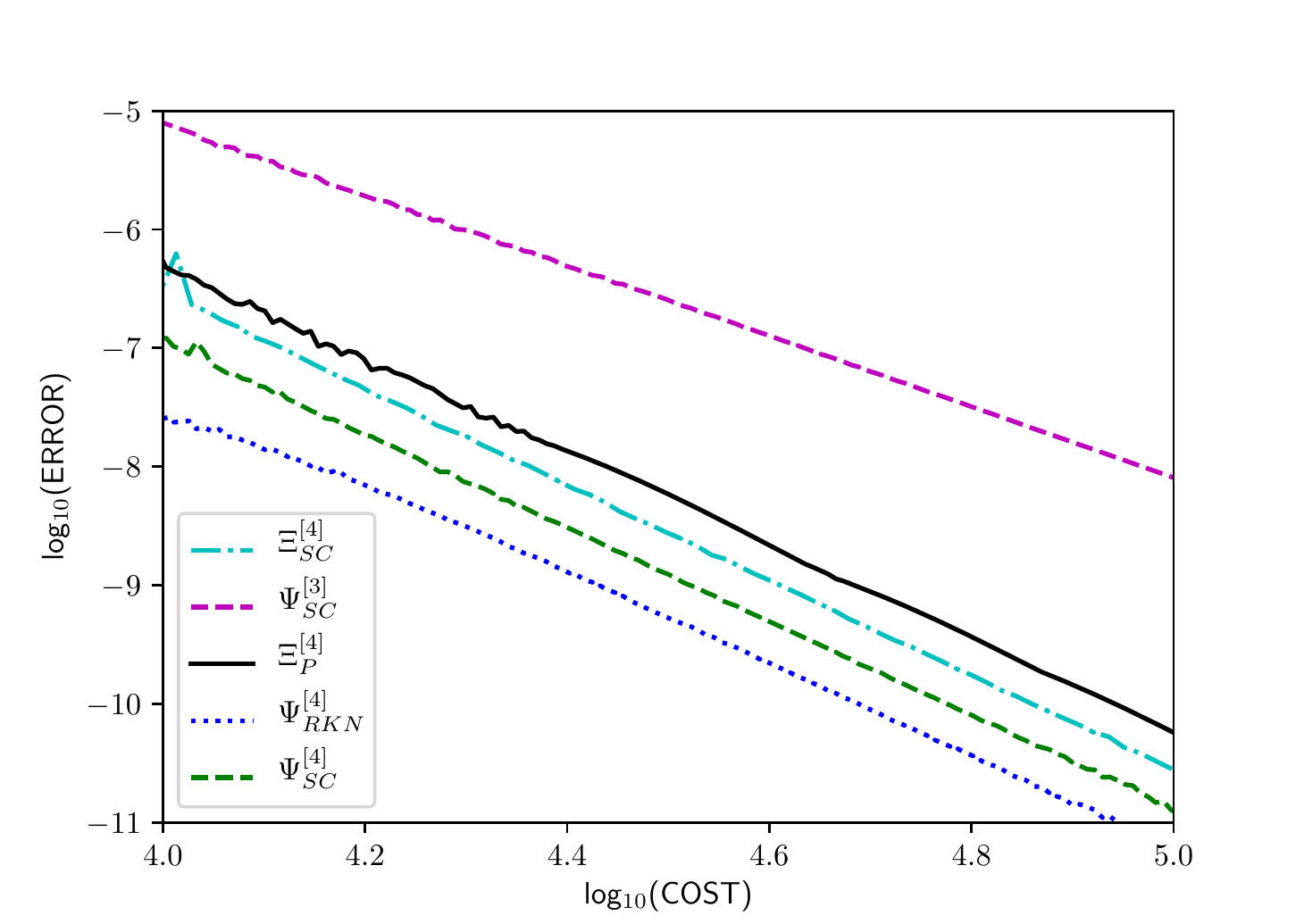}
\caption{\label{figu7} \small Maximum of error in the expected value of the energy in the interval $t \in [0, 100]$ as a function of the time step (left) and the
computational cost (number of FFTs, right) for several splitting schemes. P\"oschl--Teller potential.
}
\end{figure}

Next, we take a shorter final time $t_f = 100$ and compute the maximum error in the energy along the time interval for several step sizes $h=\Delta t$ and
integration schemes. The corresponding results are displayed in a log-log diagrama in Figure \ref{figu7} (left). The order of each method is clearly visible,
as well as the values of $h$ where instabilities take place. Finally, in Figure \ref{figu7} (right) we depict the same results but in terms of the computational
cost as measured by the number of FFTs necessary to carry out the calculations. Notice that, for this range of times, the efficiency of the 4th-order 
symmetric-conjugate methods is not far away from the optimized scheme (\ref{rkn4}) that takes into account the special property (\ref{rkn}).

\section{Concluding remarks}

Splitting and composition methods with complex coefficients have shown to be an appropriate tool in the numerical time integration of differential equations
of parabolic type, when one or more pieces of the equations are only defined in semi-groups and the aim is to get high accuracy. Since it is possible to design
methods of this class with positive real part, one is thus able to circumvent the existing order barrier for methods with real coefficients. In addition, these
methods involve smaller truncation errors than their real counterparts and also exhibit relatively large stability thresholds. On the other hand, their
computational cost notably increases, due to the use of complex arithmetic. 

More recently, it
has been shown that the particular class of symmetric-conjugate methods still exhibits remarkable preservation properties when applied to differential
equations defined by real vector fields and the solution is projected on the real axis at each integration step. Here we have extended the analysis
to problems evolving in the $\mathrm{SU}(2)$ and more generally to the numerical integration of the Schr\"odinger equation, where 
preservation of unitarity is a physical requirement. In the former case we have shown explicitly that symmetric-conjugate splitting methods are indeed
conjugate to unitary methods for sufficiently small time step sizes, and thus there is not a secular component in the unitarity error propagation.

{
With respect to the Schr\"odinger equation, the examples we collect here indicate that methods of this class (with real 
coefficients $a_j$) could safely be applied 
just as other schemes involving only real coefficients for sufficiently small step sizes, although a general theoretical analysis similar to the one developed here
for problems defined in $\mathrm{SU}(2)$ is lacking at present. Such analysis is clearly more involved, since one has to take into account
the effect of the space discretization, the possible introduction of artificial cut-off bounds for unbounded potentials, etc. In this sense, this paper
should be considered as a preliminary step for such analysis. In any case, we should remark that the use of methods with complex coefficients
in this setting does not imply any extra computational cost, since the problem has to be treated in the complex domain anyway. 
Our results show that even some of the
simplest methods within this class provide efficiencies close to the best standard splitting schemes specifically designed for the integration of the
Schr\"odinger equation. Although we have limited ourselves here to methods of order 3 and 4, it is clear that higher order integrators can also be designed,
just by solving the corresponding order conditions \cite{blanes08sac,hairer06gni}, and more efficient schemes can be obtained by taking into 
account property (\ref{rkn}) and the processing technique. It is also worth noticing that, in contrast with the time integration of parabolic differential equations,
here schemes with real and \emph{negative} coefficients $a_j$ still provide unitary approximations, and so more efficient schemes with $a_j < 0$ and
$b_j \in \mathbb{C}$ might be possible. All these issues will be treated in a forthcoming paper.
}

\subsection*{Acknowledgements}
 This work has been supported by 
Ministerio de Ciencia e Innovaci\'on (Spain) through project PID2019-104927GB-C21/AEI/10.13039/501100011033. 
A.E.-T. has been additionally funded by the predoctoral contract BES-2017-079697 (Spain).

\bibliographystyle{siam}

\end{document}